\documentclass{article}

\usepackage{amsmath}

\newcommand{\norm}[1]{\left\lVert#1\right\rVert}

\usepackage{natbib}
 \bibpunct[, ]{(}{)}{,}{a}{}{,}%
\usepackage{subcaption}
\usepackage{rotating}
\usepackage[colorlinks=true, linkcolor=blue, citecolor=blue, urlcolor=blue]{hyperref}
\usepackage{fancyvrb}
\usepackage{algorithm}
\usepackage{algorithmic}
\usepackage{tikz}
\usepackage{setspace}
\usetikzlibrary{angles, quotes} 

\newcommand{\ps}[2]{\left \langle #1, \: #2\right \rangle}
\newcommand{\abs}[1]{\left|#1\right|}
\newenvironment{Proofe}
  {\textbf{\textit{Proof.}}}
  {\hfill $ \square $}

\usepackage{tikz}
\usepackage{stmaryrd}
\newtheorem{Remark}{Remark}
\newtheorem{Lemma}{Lemma}
\newtheorem{Theorem}{Theorem}

\newtheorem{Proposition}{Proposition}
\newtheorem{Assumption}{Assumption}
\newtheorem{Definition}{Definition}
\newcommand\E{\ensuremath{{\mathbb{E}}}}
\newcommand\R{\ensuremath{{\mathbb{R}}}}
\usepackage{lipsum}
\usepackage{amssymb} 

\newtheorem{exaample}{Example}
\newenvironment{Example}
  {\begin{exaample}}
  {\hfill $ \triangleleft $\end{exaample}}

\usepackage{arxiv}

\usepackage[utf8]{inputenc} 
\usepackage[T1]{fontenc}    
\usepackage{hyperref}       
\usepackage{url}            
\usepackage{booktabs}       
\usepackage{amsfonts}       
\usepackage{nicefrac}       
\usepackage{microtype}      
\usepackage{lipsum}
\usepackage{graphicx}
\graphicspath{ {./images/} }

\title{Addressing misspecification in contextual optimization}

\author{
\begingroup\hypersetup{hidelinks}
  Omar Bennouna\thanks{Equal contribution}
  \endgroup
 \\
  Laboratory of Information and Decision Systems\\
  Massachusetts Institute of Technology\\
  \texttt{omarben@mit.edu} \\
   \And
   Jiawei Zhang\footnotemark[1]\\
  Laboratory of Information and Decision Systems\\
  Massachusetts Institute of Technology\\
  \texttt{jwzhang@mit.edu} \\
  \And
  Saurabh Amin\\
  Laboratory of Information and Decision Systems\\
  Massachusetts Institute of Technology\\
  \texttt{amins@mit.edu} \\
   \And
    Asuman Ozdaglar\\
  Laboratory of Information and Decision Systems\\
  Massachusetts Institute of Technology\\
  \texttt{asuman@mit.edu} \\
}

\begin{document}
\maketitle
\begin{abstract}
We study a linear contextual optimization problem where a decision maker has access to historical data and contextual features to learn a cost prediction model aimed at minimizing decision error. We adopt the predict-then-optimize framework for this analysis. Given that perfect model alignment with reality is often unrealistic in practice, we focus on scenarios where the chosen hypothesis set is misspecified. In this context, it remains unclear whether current contextual optimization approaches can effectively address such model misspecification. In this paper, we present a novel integrated learning and optimization approach designed to tackle model misspecification in contextual optimization. This approach offers theoretical generalizability, tractability, and optimality guarantees, along with strong practical performance. Our method involves minimizing a tractable surrogate loss that aligns with the performance value from cost vector predictions, regardless of whether the model is misspecified or not, and can be optimized in reasonable time. To our knowledge, no previous work has provided an approach with such guarantees in the context of model misspecification. 

\end{abstract}


\section{Introduction} \label{sec1}
\subsection{Problem overview}
In a number of real-world problems, the decision makers are faced with uncertainty over relevant parameters and this uncertainty can be reduced by acquiring side information. For example, consider a traffic routing problem where a fixed amount of traffic is to be allocated  over a directed graph in a way that minimizes total edge cost. The edge costs are unknown at the time of making this decision, but depend on some observable features (or context variables) that are correlated with the edge costs, such as weather, time of day, construction, and so on. Suppose that we have access to data containing previous realizations of edge costs and context variables. The decision-maker faces the task of designing a policy that maps the observed contextual information into minimum cost routing decisions. A general setting for problems of this type, typically referred to as contextual optimization problems, involves simultaneously learning a decision policy within a certain model class and solving an optimization problem using predicted costs. In this paper, we are concerned with contextual optimization with model misspecification; that is, problems in which the class of prediction models may not contain the model that recovers the true (unknown) parameters.

For simplicity, we consider that unknown parameters appear linearly in the objective function of
the stochastic optimization problem:
	\begin{equation*}
\min_{w\in W}\E_{(x,c) \sim P}(c^\top w|x)=\min_{w\in W}\E_{(x,c) \sim P}(c|x)^\top w,
	\end{equation*}
where  $c\in \mathbb{R}^d$ is a random variable correlated with the context $x\in \mathbb{R}^k$, $P$ is the joint probability distribution of $x$ and $c$ (which we will refer to as the data distribution), and $W$ is a bounded convex set. When making a decision, the decision maker does not know $c$, and only the context variable $x$ is observed. Let $\Pi$ be a set of feasible policies, i.e. a set of functions which maps a context $x\in \R^k$ to a decision $w\in \R^d$. 
  Our goal is to find a policy $\hat{\pi} \in \Pi$ approximately minimizing the cost resulting from making a decision $\hat{\pi}(x)$ given a context $x$, i.e.
\begin{equation} \label{generalcontextual}
     \hat{\pi}\in \arg \min_{\pi \in \Pi}\E_{(x,c)\sim P}(c^\top \pi(x)).
\end{equation}
This policy $\hat \pi$ involves choosing an element from the optimal solution set of problem (\ref{generalcontextual}), which therefore may not have desired continuity and smoothness properties to enable tractable parametrizations for learning an optimal policy. A commonly used approach to mitigate this issue is to “predict-then-optimize”  (\cite{hu2022fast}, \cite{elmachtoub2022smart}, \cite{sun2023maximum}). This approach, instead of learning directly a policy, involves first selecting a suitable cost vector predictor $\hat{c} : \R^k \longrightarrow \R^d$ from a given hypothesis set $\mathcal H$ (where $\hat c(x)$ is the predicted cost vector associated with context $x$),  and then using an optimization solver to determine the policy $\hat \pi$ satisfying for every context $x$, $\hat{\pi}(x) \in \arg \min_{w\in W} \hat{c}(x)^\top w $ to make decisions.

A recent active literature considers the problem of finding cost prediction models that take into account the downstream optimization problem and minimize the resulting decision error, measured as the difference between the true cost of the decision induced by the predicted cost and the optimal cost. This approach is referred to as the \emph{Integrated Learning and Optimization (ILO)} approach. Most of this literature considers the well-specified case, where the hypothesis set $\mathcal H$ contains the ground truth cost predictor $c(x)= \E(c|x)$. The misspecified case, where this ground truth is not in the hypothesis set actually arises frequently in practice  because of incomplete specification of costs or insufficient knowledge about the underlying data distribution. However, many practical situations may necessitate solving a contextual optimization problem in which
the hypothesis set does not contain the ground truth predictor. In particular, the decision maker may only
have an incomplete specification of the relationship between context and cost variables or insufficient
knowledge about the underlying data distribution. For example, in the traffic routing problem, we may not
have complete knowledge of how contextual information impacts edge costs or may not have access to
all the relevant data needed to ensure accurate cost prediction. We refer to such problems as contextual
optimization under misspecification. Despite their practical significance, these problems have not received much attention in the literature, with the notable exception of some recent papers such as  \cite{huang2024learning} and some empirical studies (see for example \cite{donti2017task}, \cite{mckenzie2023faster}, \cite{kotary2023backpropagation}), but neither of these provides global optimality properties.
\subsection{Our focus and contributions}
In this paper, we develop a computationally tractable approach to minimize the decision error (i.e., target performance) when the hypothesis set $\mathcal H$ is misspecified, using historical data. In contrast to prior work, our approach ensures global optimality. For any cost vector $\Tilde{c}\in \R^d$, we denote $w(\Tilde{c})$ a solution to the optimization problem $\min_{w\in W}\Tilde{c}^\top w$ returned by an oracle, which could be an off-the-shelf solver for example. We consider that the hypothesis set $\mathcal H$ is parametric, i.e. it can be written as $\mathcal H=\{\hat{c}_\theta,\; \theta \in \R^m\}$. We would like to minimize the following function: $ \ell_{P}(\theta)=\E_{(x,c)\sim P}(c^\top w(\hat{c}_\theta(x)))$. In general, the solution to $\min_{w\in W}\hat{c}_\theta(x)^\top w$ is not always unique, so we rewrite $\ell_P$ by accounting for the worst case solution:
\begin{equation} \label{robust formulation}
    \ell_{P}(\theta)=\E_{(x,c)\sim P}\left(\max_{w\in W^\star(\hat{c}_\theta (x))}c^\top w\right),
\end{equation}
where for every $\tilde{c}\in \R^d$, $W^\star(\tilde{c}):=\arg\min_{w\in W}\tilde{c}^\top w$. Minimizing $\ell_P$ can be quite challenging because in general it is non-convex, non-smooth, and non-continuous. Moreover, in practice, we do not have access to $P$, but rather a set of empirical samples $S:=\{(x_1,c_1),\dots,(x_n,c_n)\}$.

 Previous works (\cite{bertsimas2020predictive}, \cite{elmachtoub2022smart}, \cite{sun2023maximum}, \cite{loke2022decision}) have suggested ways to address poor regularity properties of the target performance function, leading to various procedures that aim to minimize the function but assume that $\mathcal H$ is well-specified. However, these procedures do not naturally extend to contextual optimization problems under misspecification. To the best of our knowledge, the current literature does not provide a procedure to minimize $\ell_P$ with theoretical guarantees when $\mathcal H$ is misspecified.
 
 Our work addresses these gaps in the literature by
providing a tractable procedure to minimize $\ell_P$ when $\mathcal H$ is misspecified. Our procedure enjoys both strong optimality guarantees and generalization performance. We also experimentally evaluate the performance of our procedure and show it indeed has good performance compared to mainstream methods. Importantly, our procedure leverages an intuitive surrogate loss function which naturally arises under mild assumptions on $\mathcal H$ and $P$, and does not require $\mathcal H$ to be well-specified. The proposed surrogate loss has the same set of minimizers as $\ell_P$ (see Theorem \ref{mainconsistency}). We show that our loss function generalizes well to out-of-sample distributions (i.e. optimizing the surrogate loss where $P$ is approximated by the uniform distribution over a historical dataset can yield a good solution for $\ell_P$, see Theorem \ref{surrogategen}), and it is tractable in the sense that it has no bad local minima and saddle points (see Theorem \ref{subdifferentialsurrogate}).

 Our approach contrasts with contextual optimization in the well-specified case which has been extensively studied in the literature. Previous works have mainly adopted two types of approaches: \emph{Sequential Learning and Optimization (SLO)} (\cite{hu2022fast}, \cite{bertsimas2020predictive}) and \emph{Integrated Learning and Optimization (ILO)} (\cite{elmachtoub2022smart}, \cite{donti2017task}, \cite{sun2023maximum})\footnote{We have not mentioned the Decision Rule Optimization approach since SLO and ILO perform \emph{a priori} better in theory and in practice (see \cite{elmachtoub2023estimate}).}. For both of these approaches, the general idea is to approximately optimize $\ell_P$ given a hypothesis set $\mathcal H$ and a set of past observations $S=\{(x_1,c_1),\dots,(x_n,c_n)\}$. SLO consists of focusing solely on the prediction of the cost $c$ by choosing the cost predictor with the smallest prediction error without taking into account the downstream optimization task, whereas ILO consists of finding a cost predictor which yields good decision performance rather than being a good approximation of $c$. A well-known ILO approach is to design a surrogate loss function, as notably adopted by \cite{elmachtoub2022smart}. In \cite{elmachtoub2022smart}, a new convex surrogate function is introduced (called the SPO+ loss), which is consistent with $\ell_P$ when $\mathcal H$ is well-specified, i.e. a minimizer of the expected value of their convex surrogate function,
 \begin{equation*}
     \ell_P^+(\theta):=\E_{(x,c)\sim P}\left(\ell_ {\text{SPO+}}(\hat{c}_{\theta}(x),c)\right) ,
 \end{equation*}
 where for every $c,\tilde{c}\in \R^d$,
 \begin{equation*}
	\ell_ {\text{SPO+}}(\tilde{c},c):=\max_{w\in W}(c^\top w-2(\tilde{c}^\top w-\tilde{c}^\top w(c))),
	\end{equation*}
 is also a minimizer of the target performance loss $\ell_P$. However, no theoretical optimality guarantees are provided when $\mathcal H$ is misspecified. 

 To further motivate the necessity of a new approach to address misspecification, we show in the following example that SLO and SPO+ are not suitable to handle misspecification.

\begin{Example}\label{example1}
Consider that the context variable x takes value 0 with probability 1, the underlying cost vector is $c(x)=(1,2,2)^\top \in \mathbb{R}^3$ and the hypothesis set $\mathcal H=\{\hat{c}_{\theta}(x),\;  \theta \in \mathbb R^m\}=\{c_1,c_2 \}$ where $c_1=(10,1.1,1)^\top$ and $c_2=(1,3,100)^\top$.
This misspecified problem is associated with the decision-making task modeled as a linear program:
\begin{equation*}
\min_{w_1,w_2,w_3\ge 0}c^\top w,\ \ \mathrm{s.t.}\; w_1+w_2+w_3=1.
\end{equation*}
The optimal decision for the abovementioned problem is $w(c)=(1,0,0)$. By direct calculation, we can see that the cost prediction $c_2$ yields the optimal decision $w(c)$ but $c_1$ does not. However, we have $\ell_{\text{SPO+}}(c_1,c)<\ell_{\text{SPO+}}(c_2,c)$ and $\norm{c-c_1}_2<\norm{c-c_2}_2$. Therefore, the cost prediction made by SPO+ and SLO is $c_1$, which will yield a suboptimal decision.\end{Example}

Our novel ILO approach addresses these limitations by ensuring that optimal decision performance is recovered even under model misspecification. Our surrogate loss function is directly related to an optimality condition for minimizing the target loss $\ell_P$. In fact, this surrogate is a nonnegative function and can be written as the difference between the minimum of two functions that depend on the chosen parameter $\theta$, with the property that the minimal value of the first function is equal to the minimal value of the second function if and only if $\theta$ is a minimizer of $\ell_P$. Importantly, this property holds regardless of whether $\mathcal H$ is misspecified.

 Assuming the hypothesis set is linear, our surrogate loss function has good landscape properties: every local minimum to this surrogate, is a global minimum. While it is not known whether approaches in previously mentioned works are successful in solving contextual optimization problems under misspecifiation, our approach ensures that the set of minimizers of proposed surrogate coincides with the set of minimizers of the target loss $\ell_P$ (see Theorem \ref{mainconsistency}) even when $\mathcal H$ is misspecified.
 
 Our surrogate loss function is tractable to minimize because it can be written as a difference of non-smooth convex functions (see Proposition \ref{DC structure}), we can find a stationary point using first order methods applied on its Moreau envelope. Our surrogate is tractable in the sense that every one of its stationary points is a global minimum (see Theorem \ref{subdifferentialsurrogate}). We provide theoretical guarantees bounding the gap between the value of our loss when averaging over the in-sample distribution and the out-of-sample distribution. We further give experimental evidence of our surrogate's generalization properties, and show that it outperforms SPO+ and SLO when the level of misspecification of $\mathcal H$ is reasonably high.
 \subsection{Quantiyfing the level of misspecification}
 A classical metric to quantify the level of misspecification of a given hypothesis set $\mathcal H$ is the $L^1$ or $L^2$ norm of the distance between the parameter (or the function) one seeks to estimate and the chosen hypothesis set. In our setting, it can be written as
\begin{align*}
    \gamma_{\text{class}}(\mathcal H)=\sqrt{\min_{\hat{c}\in\mathcal H}\E_{(x,c)\sim P}\left((\hat{c}(x)-c(x))^2\right)}.
\end{align*}

This metric has been previously adopted in statistics and contextual bandits (\cite{foster2020adapting}, \cite{krishnamurthy2021adapting}). We argue that $\gamma_{\text{class}}(\mathcal H)$ can be a good metric for evaluating the effect of misspecification for prediction problems, but not for contextual optimization problems in which parameter estimates enter in downstream optimization problem. This can be seen by noting that in our setting, the best achievable cost should not change when multiplying the elements of $\mathcal H$ by a positive constant, whereas $\gamma_{\text{class}}(\mathcal H)$ is not invariant to such a transformation. To address this, we introduce a new metric to quantify misspecification in contextual optimization.

\begin{Definition}\label{realizable}
We define the misspecification gap as $$\gamma_{\text{miss}}(\mathcal H)=\min_{\hat{c}\in \mathcal H}\E_{(x,c)\sim P}\left(\max_{w\in W^\star(\hat{c}(x))}c^\top w\right)-\E\left(\min_{w\in W}c^\top w\right)\geq 0.$$
The function class $\mathcal H$ is well-specified if $\gamma_{\text{miss}}(\mathcal H)=0$, and misspecified otherwise. 
\end{Definition}
Note that the first term is the smallest possible value of the target loss when choosing a cost predictor from $\mathcal H$, and the second term is the smallest possible value of the target loss in the case where $\mathcal H$ contains the ground truth cost predictor.
First, this definition of model misspecification is weaker than the previous one. Indeed, if $\gamma_{\text{miss}}(\mathcal H)>0$, then necessarily $\gamma_{\text{class}}(\mathcal H)>0$.  Second, $\gamma_{\text{miss}}(\mathcal H)$ is clearly invariant when multiplying the cost predictors by a constant. Furthermore, the misspecification gap $\gamma_{\text{miss}}(\mathcal H)$ captures the optimality gap caused by model misspecification. In other words, it is the difference between the best possible performance in the misspecified case and the best possible performance in the well-specified case.

For a given cost predictor $\hat{c}$, denoting $\gamma_{\text{miss}}(\hat{c})=\E_{(x,c)\sim P}\left(\max_{w\in W^\star(\hat{c}(x))}c^\top w\right)-\E\left(\min_{w\in W}c^\top w\right)$
and 
$\gamma_{\text{class}}(\hat{c})=\sqrt{\E_{(x,c)\sim P}\left((\hat{c}(x)-c(x))^2\right)},$
we can see that SLO focuses on bringing $\gamma_{\text{class}}( \hat{c})$ as 
close as possible to $\gamma_{\text{class}}(\mathcal H)$, whereas it is more natural to focus on bringing $\gamma_{\text{miss}}(\hat{c})$ as close to $\gamma_{\text{miss}}(\mathcal H)$ as possible, since small prediction error does not necessarily mean good decision performance in the misspecified case. In the other hand, even if current ILO methods such as SPO+ aim to bring $\gamma_{\text{miss}}(\hat{c})$ as close to $\gamma_{\text{miss}}(\mathcal H)$, it is unclear whether they can do so in the misspecified case, whereas our surrogate loss's optimality gap dominates $\gamma_{\text{miss}}(\hat c)$ (see Theorem \ref{sensitivity}).

\subsection{Previous approaches to Integrated Learning and Optimization}
A more comprehensive review of contextual optimization approaches can be found in the survey by \cite{sadana2024survey}.
\begin{paragraph}{\textbf{Directly optimizing the target loss.}}

One of the first instances of the ILO method can be traced back to \cite{donti2017task} who provide a practical way to differentiate the target loss $\ell_P$ under some regularity conditions. Similarly, others have attempted to directly minimize $\ell_P$ using some estimation of its gradient, such as unrolling (\cite{domke2012generic} \cite{monga2021algorithm}), which consists of keeping track of operations while running gradient descent in order to differentiate the final gradient descent iterate as a function of the model parameters, and implicit differentiation (\cite{amos2017optnet} \cite{agrawal2019differentiable} \cite{mckenzie2023faster} \cite{sun2022alternating}). However, $\ell_P$ is non-differentiable in general, and even if it were, there is no guarantee that it would be convex. This implies that gradient-based algorithms mentioned above do not guarantee convergence to optimum of the target loss. In contrast, our approach consists of minimizing a tractable smooth surrogate that has the same set of minimizers as $\ell_P$. This allows us to avoid the challenge coming from the lack of regularity properties for $\ell_P$. 
    
\end{paragraph}
\begin{paragraph}{\textbf{Optimizing a surrogate loss.}}
     Another increasingly popular approach, which is most relevant to ours, is optimizing the target loss using a smooth convex surrogate loss. In \cite{elmachtoub2022smart}, a new convex surrogate function, named SPO+. Minimizing SPO+ is proven to also minimize the target loss when the hypothesis set is well-specified, but no consistency results are provided when the chosen predictor class is misspecified. Furthermore, SPO+ seems to outperform SLO (i.e. yield better decision performance) when $\mathcal H$ is misspecified, since SLO only focuses on the accuracy of the prediction step, but completely disregards the performance in the optimization step, whereas ILO specifically focuses on the decision performance, and when the hypothesis set is misspecified, maximizing the prediction accuracy of the underlying cost will not necessarily result in good decision performance. Moreover, when $\mathcal H$ is well-specified, SLO methods outperform ILO (\cite{hu2022fast}, \cite{elmachtoub2022smart}). \cite{hu2022fast} give theoretical and experimental evidence showing that classical SLO methods generalize better than SPO+ when the hypothesis set is well-specified. In fact, evidence in \cite{elmachtoub2023estimate} suggests that in the well-specified case, SLO might likely have better performance than ILO approaches, and such a behavior is inverted in the misspecified case. Most previous works have not theoretically considered the misspecified case.  In \cite{huang2024learning}, a new surrogate is introduced based on an approximation of the directional derivative of $\ell_P$, and is shown to be theoretically consistent with $\ell_P$. However, only local optimality results are provided for this surrogate, whereas in our work, we provide global optimality guarantees.

     Another alternative is the surrogate introduced by \cite{sun2023maximum}, which attempts to maximize the nonbasic reduced costs (when the cost is taken equal to the predicted cost) of past realizations of ground truth optimal decisions. Such a method does not require the knowledge of historical costs, but only solutions of previously seen linear programs. In order for this surrogate to be consistent, authors assume that the chosen hypothesis set is well-specified in the sense of our definition, i.e. $\gamma_{\text{miss}}(\mathcal H)=0$. On one hand, this assumption is weaker than the usual definition of well-specification, which requires that the hypothesis set contains the ground truth predictor. On the other hand, this surrogate is only designed for linear objectives with linear constraints. Our surrogate only requires the set of feasible decisions $W$ to be convex and bounded. In a related paper, \cite{liu2021end} use a neural network structure in an inventory management problem to learn a mapping which provides the optimal merchandise order quantity and order time, and give theoretical guarantees in the well-specified setting. Other surrogates have been considered in the literature (\cite{kallus2023stochastic} \cite{loke2022decision} \cite{jeong2022exact}), which despite having good practical performance benefits, do not seem to theoretically tackle the misspecified case as opposed to our work.
\end{paragraph}

The remainder of this paper is organized as follows. In section \ref{sec2}, we formulate our problem and give the intuition behind our approach, as well as provide our main consistency theorem. We give an exact characterization of the minimizers of $\ell_P$. In section \ref{sec3}, we first establish a generalization bound of the in-sample surrogate loss $\ell_{P_n}^\beta$, where $P_n$ is the uniform distribution over the dataset $S=\{(x_1,c_1),\dots,(x_n,c_n)\}$, proving that minimizing the in-sample surrogate indeed provides a good near-optimal solution for the out-of-sample surrogate $\ell_P^\beta$, then provide a procedure to tractably minimize it, and finally show that minimizing our surrogate with our procedure indeed yields a good near-optimal solution for $\ell_P$. In section \ref{sec4}, we provide experimental evidence showing that our approach indeed performs better than state-of-the-art methods when the hypothesis set is misspecified. Finally, in section \ref{sec5}, we make some concluding remarks.

\section{Our approach}\label{sec2}

\subsection{Rewriting $\ell_P$ under a generic assumption}
Recall from (\ref{robust formulation}) that the generic form of the target loss function
in our contextual optimization problem involves the worst case value of
$c^\top w$ when $w\in \arg\min_{w\in W}\hat{c}_\theta (x)^\top w$ because the solution to $\min_{w\in W}\hat{c}_\theta (x)^\top w$ is not necessarily unique in theory. In practice, however, $\arg\min_{w\in W} \hat c_\theta(x)^\top w$ is likely to be unique. For example, in the common case where $\hat c_\theta(x)$ is a continuous random variable and $W$ is a polyhedron (or more generally when the set of directions $\Tilde{c}$ for which the solution to $\min_{w\in W}\tilde{c}^\top w$ is of Lebesgue measure equal to zero), uniqueness holds
with probability 1. Hence, due to practical reasons and theoretical
convenience, we make the following assumption henceforth:
\begin{Assumption}\label{regularity}
For any $\theta \in \mathbb R^m$ such that $\hat{c}_\theta(x) \neq 0$ almost surely, $\min_{w\in W}\hat{c}_{\theta}(x)^\top w$ has a unique solution with probability $1$ when $(x,c)\sim P$.
\end{Assumption}

We denote $W_P$ the set of measurable mappings from the support of the joint probability $P$ of $(x,c)$ to the set of feasible decisions $W$. Under Assumption \ref{regularity}, we can see that the problem of minimizing target loss can be written as follows:
\begin{align}\label{rewrite}
\min_{\theta\in \R^m}\ell_P(\theta)=&\min_{\theta\in \R^m}\min_{ w_P\in W_P}\mathbb E_{(x,c)\sim P}\left(c^\top w_P(x)\right)\\
&\text{s.t. }\forall x \in \R^k,\; w_P(x)\in \arg \min_{w\in W}\hat{c}_\theta(x)^\top w.
\end{align}

In the formulation above, the minimum was taken outside of the expectation because, to minimize the function $ w_P \longmapsto \E_{(x,c)\sim P}\left(c^\top w_P(x)\right) $ for $ w_P \in W_P$, it is sufficient to select $ w_P(x) $ as the minimizer of $ w \longmapsto c^\top w $ for every possible realization of $x$. The following proposition provides a sufficient condition for Assumption \ref{regularity} to hold when $W$ is a polyhedron.
\begin{Proposition} \label{continuousdist}
If the context variable x has a continuous probability distribution and the mapping $x\longmapsto \hat{c}_\theta (x)$ is a nonzero analytic function for any $\theta\neq 0$, then Assumption \ref{regularity} holds for any polyhedron $W$ that satisfies $\hat{c}_{\theta}(x)^\top v\neq 0$ for any $\theta\in \R^m \setminus\{0\}$ and any nonzero $v\in W$.
\end{Proposition}
\begin{Proofe}
Let $\{S_i\}_{i\in I}$ be all the faces of $W$, where $I$ is a finite index set. Let $v_i$ be an arbitrary tangent direction of $S_i$ and let $\mathcal{V}=\{v_i\}_{i\in I}$ be a finite set.
Let $U$ be the set containing all $u$ such that the solution to the linear program $\min_{w\in W}u^\top w, u\in U$ is not unique.
Then any $u\in U$ satisfies $u^\top v_i=0$ for some $i$.
Therefore, $\hat{c}_{\theta}(x)\in U$ if $\hat{c}_{\theta}(x)^\top v_i=0$ for some $i$.
Let $X_0:=\{x\in \R^k\mid \hat{c}_{\theta}(x)\in U\}$ and $X_i:=\{x\in \R^k\mid \hat{c}_{\theta}(x)^\top v_i=0\}$. We have $X_0\subseteq \bigcup_{i\in I}^{}X_i$.
By the assumption in the Proposition, $\hat{c}_{\theta}(x)^\top v_i$ is a nonzero analytic function and then by \cite{mityagin2015zero}, $X_i$ is a zero-measure set.
Therefore, $X_0=\bigcup_{i\in I}^{}X_i$ is a zero-measure set since $I$ is finite.
Finally, since $X_0$ is of measure zero and $x$ is a continuous random variable, $x$ lies outside of $X_0$ almost surely, i.e. the solution to $\min_{w\in W}\hat{c}_\theta(x)^\top w$ is unique almost surely.
\end{Proofe}

Our task now becomes to solve the bilevel optimization problem (\ref{rewrite}). Since known approaches to solving bilevel problems cannot be directly applied to solve (\ref{rewrite}), we develop an approach based on a surrogate loss function that enjoys two properties:

\begin{enumerate}
\item \textbf{Consistency:} The optimal solutions for the surrogate loss are also optimal for the target loss;
\item \textbf{Tractability:} The surrogate loss is tractable to optimize.
\end{enumerate}
\subsection{Introducing our surrogate loss function}
We denote $\beta^\star _{\mathcal H,P}=\min_{\theta \in \R^m}\ell_P(\theta)$, $\beta_{\text{max},P}=\E_{(x,c)\sim P}(\max_{w\in W}c^\top w)$, and $\beta_{\text{min},P}=\E_{(x,c)\sim P}\left(\min_{w\in W}c^\top w\right)$. Let $\beta \in \R$. We first study the property $\ell_P(\theta)\le \beta$ when $\beta$ satisfies $\beta_{\mathcal H,P}^\star \leq \beta < \beta_{\text{max},P}$ and $\hat{c}_\theta(x)\neq 0$ almost surely. The inequality $\beta^\star_{\mathcal H,P}\leq \beta$ ensures that the condition $\ell_P(\theta)\leq \beta$ is feasible, and the inequality $\beta < \beta_{\text{max},P}$ ensures that the condition $\ell_P(\theta)\leq \beta$ is not trivial. Notice that if we take $\beta=\beta^\star_{\mathcal H,P}$, any $\theta \in \R^m$ satisfying $\ell_P(\theta)\leq \beta$ is a minimizer of $\ell_P(\theta)$. When Assumption \ref{regularity} holds, and when $\hat{c}_\theta(x)\neq 0$ almost surely, $W^\star(\hat{c}_{\theta}(x))$ contains only one element for almost every $x$. 
Hence, we denote $w(\hat{c}_{\theta}(x))$ to be the unique solution in $W^\star(\hat{c}_{\theta}(x))$ when $W^\star(\hat{c}_{\theta}(x))$ contains only one element.
The condition $\ell_P(\theta)\le \beta$ can be rewritten as $\mathbb{E}_{(x,c)\sim P}(c^\top w(\hat{c}_{\theta}(x)))\le \beta$ when $\hat{c}_\theta (x)\neq 0$ almost surely. For a given $\theta\in \R^m$ and $w_P\in W_P$, if $w_P$ satisfies
\begin{equation}\label{without}
w_P\in \arg\min_{w_P\in W_P}\mathbb{E}_{(x,c)\sim P}(\hat{c}_{\theta}(x)^\top  w_P(x))
,
\end{equation}
then we have, $w_P(x) = w(\hat{c}_{\theta}(x))$ almost surely, as for each $x$, the unique minimizer of $\hat{c}_{\theta}(x)^\top w$ over $w \in W$ is given by $w\left(\hat{c}_\theta(x)\right)$, according to Assumption \ref{regularity}.

If $\theta\in \R^m$ satisfies $\ell_P(\theta)\leq \beta$, and $w_P$ is the measurable mapping satisfying condition (\ref{without}) (and hence $w_P(x)=w(\hat{c}_{\theta}(x))$ almost surely), then we have $\mathbb E_{(x,c)\sim P}\left(c^\top w_P\left(x\right)\right)   \leq \beta$. This suggests that adding the linear constraint $\mathbb E_{(x,c)\sim P}\left(c^\top w_P\left(x\right)\right)   \leq \beta$ to  \eqref{without} does not change the set of minimizers of $w_P\longmapsto \mathbb{E}_{(x,c)\sim P}(\hat{c}_{\theta}(x)^\top  w_P(x))$ when $\ell_P(\theta)\leq \beta$.
Defining $\overline{W}_P^\beta=\left\{\overline{w}_P\in W_P ,\; \mathbb E_{(x,c)\sim P}\left(c^\top  \overline{w}_P(x)\right) \leq \beta \right\}$, our latter statement means that if $\ell_P(\theta)\leq \beta$, the two optimization problems
\begin{equation}\label{without2} \min_{w_P\in W_P}\mathbb{E}_{(x,c)\sim P}\left(\hat{c}_{\theta}(x)^\top w_P(x)\right)
\end{equation}
and
\begin{equation}\label{with2}
\min_{\overline w_P\in \overline W_P^\beta}\mathbb{E}_{(x,c)\sim P}\left(\hat{c}_{\theta}(x)^\top \overline w_P(x)\right)
\end{equation}
would result in identical values of the objective function. The above analysis suggests that a possible approach is to use the difference of the optimal function values of the two above optimization problems to be the surrogate function. We introduce a new surrogate loss function -- referred as Consistent Integrated Learning and Optimization (CILO) loss -- for contextual optimization under misspecification.
\begin{Definition}{(CILO loss)} \label{CILO def}
For $\beta \in \R$, we define the function $\ell_P^\beta$ as
\begin{equation*}
\forall \theta \in \R^m, \;\ell_P^\beta(\theta):=\min_{\overline{w}_P^{\beta}\in \overline{W}_P^{\beta}}\mathbb E_{(x,c)\sim P}\left(\hat{c}_\theta (x)^\top \overline{w}_P^{\beta}(x)\right) - \min_{w_P\in W_P}\mathbb E_{(x,c)\sim P}\left(\hat{c}_\theta (x)^\top w_P(x)\right).
\end{equation*}
\end{Definition}
The following example shows that whereas SLO and SPO+ fail to provide a minimizer for $\ell_P$, this new surrogate successfully does so.

\begin{Example}\label{example2}
We revisit example \ref{example1}. We have seen that SPO+ and SLO do not pick the best element of the hypothesis set $\mathcal H$. If we choose $\beta=1$. For a given $\theta\in \R^m$, we denote $\hat{c}_\theta$ as the value of $\hat{c}_\theta(x)$ (which does not depend on $x$) for any $x$.
\begin{align}
    \ell_P ^\beta(\theta)=&\min_{} \hat{c}^\top_\theta w - \hat{c}^\top_\theta w(\hat{c}_\theta)\\ 
    \text{s.t.} &\;w_1,w_2,w_3\geq 0\\
    &\;w_1+w_2+w_3=1 \label{linear constrains}
    \\&\; w_1+2w_2+2w_3=1 \label{constraint beta}
\end{align}
Here, (\ref{linear constrains}) and (\ref{constraint beta})  are the constraints involved in the definition of $\overline W_P^\beta$ (which in the definition are written as $\forall x \in \R^k,\; w_P(x)\in W$ and $\E_{(x,c)\sim P}(c^\top w_P(x))=\beta$). We denote $\theta_1,\theta_2\in \R^m$ such that $\hat{c}_{\theta_1}=c_1$ and $\hat{c}_{\theta_2}=c_2$.
\begin{align*}
    \ell^{\beta }_P(\theta_1)&=c_1^\top (1,0,0) - c_1^\top (0,0,1)=9\\
    \ell^{\beta}_P(\theta_2)&=c_2^\top (1,0,0) - c_2^\top (1,0,0)=0<\ell^{\beta }_P(\theta_1).
\end{align*}
We see that indeed here, the global minimum of this new loss is optimal for the target loss. 
\end{Example}

We will now prove that the consistency observed in example \ref{example2} holds in general. The lemma below will be useful when proving our main consistency theorem. \begin{Lemma} \label{positivity}
For every $\beta \in \R$ such that $\beta^\star_{\mathcal H,P} \leq \beta$, $\ell_P^\beta$ is non-negative.
\end{Lemma}
By the above analysis, we know that any $\theta\in \R^m$ such that $\ell_P(\theta)\le \beta$ satisfies $\ell_P^{\beta}(\theta)=0$.
Interestingly, we can show the reverse direction: if $\ell_P^{\beta}(\theta)=0$ and $\hat{c}_{\theta}(x)$ is nonzero almost surely, then $\ell_P(\theta)\le \beta$.
\begin{Theorem} \label{mainconsistency}
    Let $\beta \in \R$ such that $\beta^\star_{\mathcal H,P} \leq \beta < \beta_{\text{max},P}$. Under Assumption \ref{regularity}, for every $\theta \in \R^m$ with $\hat{c}_{\theta}(x)\ne 0$ almost surely, $\ell_P(\theta)\leq \beta$ if and only if $\theta$ is a minimizer of $\ell_P^\beta$. In particular, when $\beta =\beta^\star_{\mathcal H,P}$, $\theta$ is a minimizer of $\ell_P^\beta$ if and only if it is a minimizer of $\ell_P$. 
\end{Theorem}
\begin{Proofe}
    Let $\beta \in \R$ such that $\beta^\star_{\mathcal H,P} \leq \beta < \beta_{\text{max},P}$ and $\theta \in \R^m$ satisfying $\hat{c}_{\theta}(x)\ne 0$ almost surely. 
    Assume that $\theta$ is optimal for $\ell_P^\beta$, i.e. $\ell_P^\beta(\theta)=0$.
Let $\overline w_P(x)$ be the solution of \eqref{with2}. Then $\ell_P^{\beta}(\theta)=0$ implies
\begin{equation*}\label{intermediatecondition}
    \E_{(x,c)\sim P}\left(\hat{c}_\theta(x)^\top \overline w_P(x) - \hat{c}_\theta(x)^\top w\left(\hat{c}_\theta(x)\right) \right)=0.
\end{equation*}
Since the random variable $\hat{c}_\theta(x)^\top \overline w_P(x) - \hat{c}_\theta(x)^\top w\left(\hat{c}_\theta(x)\right)$ is positive almost surely, the equality above implies that $\hat{c}_\theta(x)^\top \overline w_P(x) - \hat{c}_\theta(x)^\top w\left(\hat{c}_\theta(x)\right)=0$ almost surely, i.e. $\overline w_P$ is also a minimizer of $w_P \longmapsto  E_{(x,c)\sim P}\left(\hat{c}_\theta (x)^\top w_P(x)\right)$, and hence under Assumption \ref{regularity}, we have $\overline w_P(x)=w(\hat{c}_{\theta}(x))$ almost surely. This finally gives 
\begin{align*}
    \ell_P(\theta)= \E_{(x,c)\sim P}(c^\top w\left(\hat{c}_\theta(x)\right)) =\E_{(x,c)\sim P}(c^\top \overline w_P(x))\leq \beta.
\end{align*}

The other direction is already shown before during the development of $\ell_P^{\beta}$.
\end{Proofe}

This result also suggests that a natural approach to minimizing $\ell_P$ is to first find a tractable method for minimizing the CILO loss $\ell_P^\beta$ for an appropriate choice of auxiliary parameter $\beta$, while ensuring that the candidate minimizer $\theta$ satisfies $\hat{c}_\theta(x) \neq 0$ almost surely. Although $\beta^\star_{\mathcal{H}, P}$ is unknown, one might expect that if the optimization procedure can achieve a sufficiently small value of $\ell_P^\beta$, the corresponding minimizer would likely ensure that the original target loss function $\ell_P$ is less than $\beta$. To select the auxiliary parameter $\beta$, one can perform a line search within a suitable interval $[\underline{\beta}, \overline{\beta}]$ and choose the parameter that yields the smallest value of the CILO loss. We will now proceed to develop a tractable procedure to optimize $\ell_P^\beta$ (so that the line search for the best $\beta$ is also feasible) and formalize the relationship between $\ell_P^\beta$ near optimality and the corresponding value of the target loss $\ell_P$. In doing so, we address the key limitation of current literature in contextual optimization -- namely, previously proposed surrogate loss functions, while tractable,
only satisfy theoretical optimality under well-specified settings. These approaches may provide reasonably good experimental performance for some misspecified settings; however, there currently does not exist an approach that is guaranteed to ensure
global optimality for \emph{any} level of misspecification. We bridge this gap in Section \ref{sec3}.

\section{Technical approach} \label{sec3}
In this section, we formalize our technical approach based on the CILO loss function and the consistency result. We first mention three key issues that our approach seeks to address. In the remainder of the paper, we denote by $\norm{.}$ the $L^2$ norm.

Firstly, we do not have access to the joint distribution P of random
variables $(x,c)$; instead we assume that we have access to a historical
dataset $S=\{(x_1,c_1),\dots,(x_n,c_n)\}$, where $n$ is the number of samples and each sample $(x_i,c_i)$, $i\in [n]$ is sampled from $P$. Thus, we seek to optimize the
empirical version of CILO, denoted $\ell_{P_n}^\beta$, where $P_n$ is the uniform distribution over the dataset $S$. The question we need to address is whether by minimizing empirical CILO, we can obtain good generalization performance (i.e.,
guarantees on out-of-sample CILO loss). Theorem \ref{surrogategen} shows
that this is indeed the case.

Secondly, the CILO loss $\ell_P^\beta$ is a non-convex and non-smooth function, and we need a technique to optimize it in a tractable manner. We address this issue by leveraging Moreau envelope smoothing technique (\cite{sun2021algorithms}) to transform the problem of minimizing $\ell_P^\beta$ into another optimization problem whose
objective function enjoys good landscape; in particular it is smooth and has no ``bad'' first-order stationary points or local minima (see Theorem \ref{subdifferentialsurrogate}). Consequently, this smoothed optimization problem is conducive to gradient descent.

Thirdly, to guarantee consistency, we must ensure that the optimization procedure is able to find minimizer $\theta$ that verifies $\hat c_\theta(x)\neq$ 0 almost surely (see Theorem \ref{mainconsistency}). We address this issue by refining our smoothing procedure so that gradient descent on smoothed CILO loss results in a minimizer that verifies $\hat c_\theta(x)\neq 0$ almost surely (see Theorem \ref{move}).

These three steps together ensure that we have a tractable approach to minimizing empirical CILO loss, resulting in solution that has good generalization perform
optimality guarantee.

\subsection{Generalization performance} \label{sec gen}
We study the generalization performance of the empirical version of CILO loss and show that by optimizing the empirical CILO loss $\ell_{P_n}^\beta$, we can ensure a nearly optimal value for its out of sample counterpart $\ell_P^\beta$.
We first make the following boundedness assumptions.
\begin{Assumption} \label{boundedness}
$W$ is closed and bounded. We denote $B_W=\sup_{w\in W}\norm{w}$.
\end{Assumption}
\begin{Assumption} \label{boundedness2}
    There exists $K\geq0$ such that for all $x\in \R^k$, $\norm{c(x)}\leq K$.
\end{Assumption}
\begin{Assumption} \label{normboundedness}
    When $\theta=0$, $\hat{c}_\theta(x)=0$. Furthermore, there exists $B_{\Phi}\geq 0$ such that the gradient with respect to $\theta$, $\nabla \hat{c}_\theta (x)$, is piecewise continuous and bounded by $B_\Phi$ for every $\theta \in \R^m$ and $x \in \R^k$.
\end{Assumption}
The first part of Assumption \ref{normboundedness} is natural when we assume that there exists $\theta_0 \in \R^m$ such that $\hat{c}_{\theta_0}=0$. By reparametrizing the hypothesis set as $\mathcal{H}'={\hat{c}_{\theta+\theta_0}, \theta \in \R^m}$, we obtain $\hat{c}_0=0$. Furthermore, a direct consequence of this assumption is that for every $x \in \R^k$ and $\theta \in \R^m$,
\begin{align*}
    \norm{\hat{c}_\theta (x)} = \norm{\int_{0}^{1}\nabla \hat{c}_{t\theta}(x)^\top \theta  dt}\leq \int_{0}^{1}\norm{\nabla \hat{c}_{t\theta}(x)^\top \theta}  dt\leq B_\Phi \norm{\theta}.
\end{align*}

We now present our main generalization result.
\begin{Theorem} \label{surrogategen}
    Let $\beta \geq \beta^\star_{\mathcal H,P}$. Let $\theta^\star \in \R^m$ such that $\ell_{P_n}^\beta(\theta^\star)\leq \varepsilon$. Assume that there exists $D\geq 0$ such that $\norm{\theta^\star} \leq D$. Under assumptions \ref{boundedness}, \ref{boundedness2}, and \ref{normboundedness}, and assuming $\beta^\star_{\mathcal H,P}>\beta_{\text{min},P}$, with probability at least $1-\delta$, we have 
    \begin{align*}
        \ell_{P}^{\beta } (\theta^\star)\leq \varepsilon + O\left(\frac{1}{\gamma_{\text{miss}}(\mathcal H) + \beta -\beta_{\mathcal H,P}^\star} \sqrt{\frac{\log \frac{1}{\delta}}{n}}\right).
    \end{align*}
\end{Theorem}
 Notice here that when the hypothesis set $\mathcal H$ is misspecified, we have $\gamma_{\text{miss}}(\mathcal H)>0$, and consequently the term $\frac{1}{\gamma_{\text{miss}}(\mathcal H) + \beta -\beta_{\mathcal{H},P}^\star}$ is bounded. If $\mathcal{H}$ is (nearly) well-specified, we can choose a larger $\beta$ to guarantee the generalization bound to be good.
This theorem implies that our surrogate loss can generalize when $n$ is large and hence minimizing the empirical version of the surrogate loss $\ell_{P_n}^\beta$ yields small optimality gap for its out-of-sample counterpart $\ell_P^{\beta}$.  
Using Theorem \ref{mainconsistency} and Theorem \ref{surrogategen}, we can deduce that finding $\theta$ such that $\hat{c}_{\theta}(x)\ne 0$ almost surely and minimizing $\ell_{P_n}^{\beta}$ yields an optimal value for $\ell_P(\theta)$.

Theorem \ref{mainconsistency} ensures that the closer $\theta$ is to optimality for $\ell_P^\beta$ for a well-chosen $\beta \in \R$, the closer it is to optimality for $\ell_P$. When $\theta$ is nearly optimal for $\ell_P^\beta$, we aim to quantify how close it is to optimality for $\ell_P$. By adopting a stability assumption regarding the linear optimization problem, we can establish a more concrete accuracy bound for the solution obtained from the empirical version $\ell_{P_n}^{\beta}$.

\subsection{Relationship between the optimality gaps of $\ell_P$ and $\ell_P^\beta$}
For ease of presentation,  we make the following assumption from now on:
\begin{Assumption}\label{characterize zero}
For $\theta\ne 0$, $\hat{c}_{\theta}(x)\ne 0$ almost surely.
\end{Assumption}
Assumption \ref{characterize zero} holds when for every $\theta\in \R^m\setminus\{0\}$, $\hat{c}_\theta$ is a nonzero analytic function and $x$ has a continuous distribution. In this case, the set of zeroes of $\hat{c}_\theta$ is of measure equal to zero (see \cite{mityagin2015zero}) and consequently $x$ does not belong to this set almost surely. This assumption, coupled with Assumption \ref{regularity}, implies that for any $\theta \in \R^m\setminus\{0\}$, the problem $\min_{w\in W} \hat{c}_\theta(x)^\top w$ has a unique solution almost surely.

\begin{figure}[h]
    \centering
    \begin{subfigure}[b]{0.45\textwidth}
        \centering
        \begin{tikzpicture}[scale=1.65]
    \coordinate (A) at (0,0);
    \coordinate (B) at (2,0);
    \coordinate (C) at (3,1.5);
    \coordinate (D) at (2,3);
    \coordinate (E) at (0,3);
    \coordinate (F) at (-1,1.5);

    \draw[black] (A) -- (B) -- (C) -- (D) -- (E) -- (F) -- cycle;

    \node at (1,1.3) [text=blue] {$\hat{c}_\theta(x)$};
    \node at (0,0.4) [text=blue] {\tiny $\omega(\theta,x)$};
    \node at (A) [text=blue,below] {\small $w(\hat{c}_\theta(x))$};
    \node at (F) [text=blue,left] {\small $\underline{w}(\theta,x)$};
    \node at (-1,0.5) {\huge $W$};
    \draw[blue, ->] (A) -- (1,2);
    \draw[dashed] (A) -- (2,2/1.5);
    \draw[red, dashed] (A) -- (1.7,1.3);
    \draw[red,dashed] (A) -- (1.854,1.069);
    \fill[red, opacity=0.3] (A) -- (1.854,1.069) -- (1.7,1.3) -- cycle;
    \draw[dashed] (F) -- (0.999, 0.166);
    \draw[red, dashed] (F) -- (0.699, 0.200);
        \draw[red,dashed] (F) -- (0.853, 0.431);
    \fill[red, opacity=0.3] (F) -- (0.69946154, 0.20069231) -- (0.85346154, 0.43169231) -- cycle;
\draw[dashed] (E) -- (-0,0.596);
\draw[dashed,red] (E) -- (-0.139,0.865);
\draw[dashed,red] (E) -- (0.138,0.865);
\fill[red, opacity=0.3] (E) -- (-0.139,0.865)-- (0.138,0.865) -- cycle;
\draw[dashed] (D) -- (0,1.667);
\draw[dashed,red] (D) -- (0.146,1.931);
\draw[dashed,red] (D) -- (0.3,1.7);
\fill[red, opacity=0.3] (D) -- (0.146,1.931)-- (0.3,1.7) -- cycle;
\draw[dashed] (C) -- (1,2.834);
\draw[dashed,red] (C) -- (1.32,2.8);
\draw[dashed,red] (C) -- (1.147,2.569);
\fill[red, opacity=0.3] (C) -- (1.32,2.8)-- (1.147,2.569) -- cycle;
\draw[dashed] (B) -- (2,2.4);
\draw[dashed,red] (B) -- (2.139,2.134);
\draw[dashed,red] (B) -- (1.861,2.135);
\fill[red, opacity=0.3] (B) -- (2.139,2.134)-- (1.861,2.135) -- cycle;
\coordinate (P) at (1,2);  
    \draw pic[ draw=blue, angle eccentricity=1.2, angle radius=0.5cm] {angle=P--A--F};

\end{tikzpicture}
        \caption{Less likely directions for $\hat{c}_\theta(x)$ relative to $W$ when it is a polyhedron in $\R^2$}
        \label{fig:1a}
    \end{subfigure}
    \hfill
    \begin{subfigure}[b]{0.45\textwidth}
        \centering
        \begin{tikzpicture}[scale=1.2]
\draw[dashed] (1.0,1.5) -- (3.0,2.833);
\draw[dashed,red] (1.0,1.5) -- (2.7,2.8);
\draw[dashed,red] (1.0,1.5) -- (2.854,2.569);
\draw[dashed] (1.0,1.5) -- (1.0,3.9);
\draw[dashed,red] (1.0,1.5) -- (1.1389999999999998,3.634);
\draw[dashed,red] (1.0,1.5) -- (0.861,3.635);
\draw[dashed] (1.0,1.5) -- (-1.0,2.834);
\draw[dashed,red] (1.0,1.5) -- (-0.6799999999999999,2.8);
\draw[dashed,red] (1.0,1.5) -- (-0.853,2.569);
\draw[dashed] (1.0,1.5) -- (-1.0,0.16700000000000004);
\draw[dashed,red] (1.0,1.5) -- (-0.854,0.43100000000000005);
\draw[dashed,red] (1.0,1.5) -- (-0.7,0.19999999999999996);
\draw[dashed] (1.0,1.5) -- (1.0,-0.904);
\draw[dashed,red] (1.0,1.5) -- (0.861,-0.635);
\draw[dashed,red] (1.0,1.5) -- (1.138,-0.635);
\draw[dashed] (1.0,1.5) -- (2.999,0.166);
\draw[dashed,red] (1.0,1.5) -- (2.699,0.2);
\draw[dashed,red] (1.0,1.5) -- (2.8529999999999998,0.431);
\fill[red, opacity=0.3] (1.0,1.5) -- (2.854,2.569)-- (2.7,2.8) -- cycle;
\fill[red, opacity=0.3] (1.0,1.5) -- (0.861,3.635)-- (1.1389999999999998,3.634) -- cycle;
\fill[red, opacity=0.3] (1.0,1.5) -- (-0.853,2.569)-- (-0.6799999999999999,2.8) -- cycle;
\fill[red, opacity=0.3] (1.0,1.5) -- (-0.7,0.19999999999999996)-- (-0.854,0.43100000000000005) -- cycle;
\fill[red, opacity=0.3] (1.0,1.5) -- (1.138,-0.635)-- (0.861,-0.635) -- cycle;
\fill[red, opacity=0.3] (1.0,1.5) -- (2.8529999999999998,0.431)-- (2.699,0.2) -- cycle;
    
\end{tikzpicture}
        \caption{Less likely directions for $\hat{c}_\theta(x)$}
        \label{fig:1b}
    \end{subfigure}
    \caption{Intuition for assumption \ref{marginal hypothesis}: $\hat{c}_\theta(x) $ is less likely to have a direction included in the red cones, which represent nearly perpendicular directions to one of the faces of the polyhedron.}
    \label{fig:main}
\end{figure}
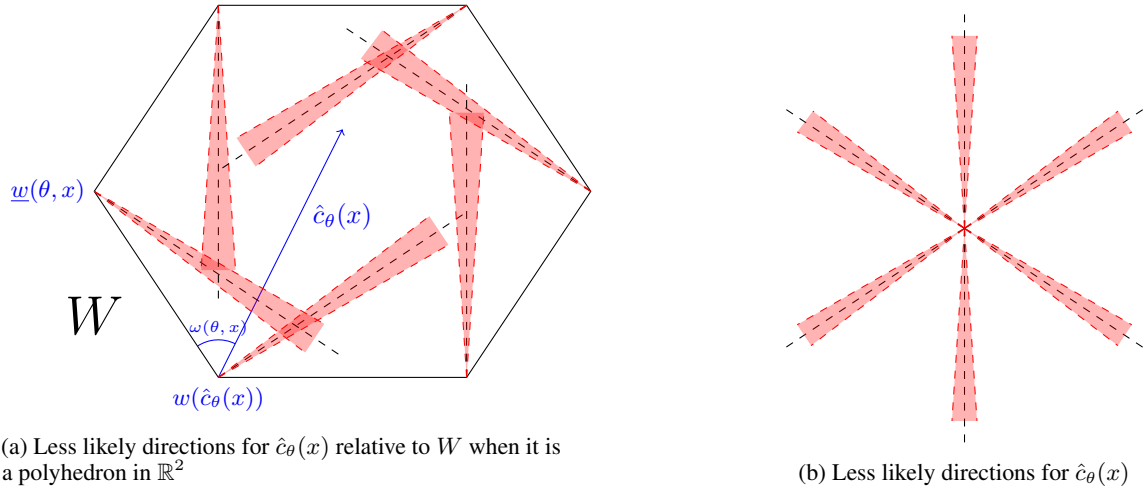

We give a more concrete optimality guarantee obtained by minimizing $\ell_{P}^{\beta}$. We first mention the assumption in \cite{hu2022fast} which enables the authors to get good generalization guarantees for SLO in the well-specified case. Note that we will not be making this assumption, but we mention it just for reference.
\begin{Assumption} \label{kallusregularity}
    Assume that $W$ is a polyhedron. We denote $W^\angle$ the set of extreme points of $W$. For a given context $x\in \R^k$, we denote 
    \begin{align*}
        \Delta(x)=\begin{cases}
            \min_{w \in W^\angle \setminus W^\star (c(x))}c(x)^\top w - \min_{w \in W^\angle}c(x)^\top w  & \text{ if }W^\star (c(x))\neq W\\
            0  & \text{ else. }
        \end{cases}
    \end{align*}
    Assume that for some $\alpha,\gamma\geq 0$,
    \begin{align*}
        \forall t>0,\;\mathbb P(0<\Delta(x)\leq t)\leq \left(\frac{\gamma t}{B_W}\right)^\alpha.
    \end{align*}
\end{Assumption}
 We make a similar assumption, which is independent of the ground truth predictor, but rather depends on the chosen hypothesis set.
\begin{Assumption} \label{marginal hypothesis}
    Assume that $W$ is a polyhedron. We denote $W^\angle$ the set of extreme points of $W$. For a given context $x\in \R^k$, we denote for every $\theta \in \R^m$, 
    \begin{align*}
        \Delta_\theta(x)=\begin{cases}
             \min_{w \in W^\angle \setminus W^\star (\hat{c}_\theta(x))}\hat{c}_\theta(x)^\top w- \min_{w \in W^\angle}\hat{c}_\theta(x)^\top w  & \text{ if }W^\star (\hat{c}_\theta(x))\neq W\\
            0  & \text{ else. }
        \end{cases}
    \end{align*}
    Assume that for some $\alpha>0$ and $\gamma\geq 0$, for every $\theta \in \R^m$
    \begin{align*}
        \forall t>0,\;\mathbb P(0<\Delta_\theta(x)\leq \norm{\theta}t)\leq \left(\frac{\gamma t}{B_W}\right)^\alpha.
    \end{align*}
\end{Assumption}
The assumption above ensures the stability of the target loss $\ell_P$. Specifically, it strengthens the uniqueness guarantee provided by Assumption \ref{regularity}. Besides ensuring that $\min_{w \in W} \hat{c}_\theta(x)^\top w$ has a unique solution almost surely, it offers a sharper measure for how sensitive the mapping $w \longmapsto \hat{c}_\theta(x)^\top w$ is to deviations from $w(\hat{c}_\theta(x))$ with high probability. This is a reasonable assumption when the distribution of $\hat{c}_\theta(x)$ is continuous. Indeed, when $W$ is a polyhedron, denoting for $\theta \neq 0$, $\underline w(\theta,x)=\arg\min_{w\in W^\angle \setminus\{w(\hat{c}_\theta(x))\}}\hat{c}_\theta(x)^\top w$, and
\begin{align*}
    \cos \left(\omega(\theta,x)\right) = \frac{\hat{c}_\theta(x)^\top (\underline w(\theta,x) - w(\hat{c}_\theta(x))}{\norm{\hat{c}_\theta(x)}\norm{ \underline w(\theta,x) - w(\hat{c}_\theta(x))}},
\end{align*}
 Assumption \ref{marginal hypothesis} can be rewritten as 
 \begin{align} \label{equivalent ineq}
     \mathbb P\left(\norm{\hat{c}_\theta(x)}\norm{\underline w(\theta,x) -  w(\hat{c}_\theta(x))} \cos(\omega(\theta,x))\leq \norm{\theta}t\right)\leq \left(\frac{\gamma t}{B_W}\right)^\alpha.
 \end{align}
We do not need to focus on the term $\norm{ \underline w(\theta,x)-w(\hat{c}_\theta(x))}$, as it is both upper bounded and bounded away from zero. Therefore, inequality (\ref{equivalent ineq}) is equivalent to $\hat{c}_\theta(x)$ having a norm that is bounded away from zero with high probability, and an angle $\omega (\theta,x)$ that is bounded away from $\frac{\pi}{2}$ when $\theta$ is bounded away from $0$. In other words, $\hat{c}_\theta(x)$ is likely to have a direction that is not too close to being perpendicular to one of the faces of the polyhedron $W$--that is, the probability of this direction to be falling within one of the red cones shown in Figure \ref{fig:main} decays to $0$ as the cones get more narrow. When $\hat{c}_\theta(x)$ follows a continuous distribution, such a property on the direction of $\hat{c}_\theta(x)$ is reasonable.

In contrast with our assumption, we have no control on whether Assumption \ref{kallusregularity} is verified or not, and also cannot verify if it is satisfied unless we have access to the ground truth predictor $c(x)$. Our assumption depends only on the chosen hypothesis set, which is something the decision maker can control. A consequence of this assumption is an inequality providing a more precise quantification of how optimality transferred from $\ell_{P}^\beta$ to $\ell_P$.


The following theorem provides a relationship between optimality gaps of $\ell_P^{\beta}(\theta)$ and $\ell_P(\theta)$ when $\theta$ is bounded away from $0$, i.e. there exists $u>0$ such that $\norm{\theta}\geq u$.
\begin{Theorem} \label{sensitivity}
    Under assumptions \ref{boundedness}, \ref{boundedness2}, \ref{normboundedness} and \ref{marginal hypothesis}, there exists $\alpha>0$ such that for all $\theta \in \R^m$ bounded away from $0$ and $\beta \geq \beta_{\mathcal H, P}^\star$,
    \begin{align*}
         \ell_P(\theta) \leq \beta + O\left(\ell_P^\beta(\theta)^{1-\frac{1}{1+\alpha}}\right).
    \end{align*}

\end{Theorem}
The proof of this theorem can be found in appendix \ref{sensitivity proof}. Assuming that $\theta$ is bounded away from 0 is crucial here, because when $\theta = 0$, the uniqueness property in Assumption \ref{regularity} is not guaranteed, and the consistency in Theorem \ref{mainconsistency} is not satisfied. Similarly, when $\theta$ is close to 0, $\ell_P^\beta$ exhibits poor behavior, as it approaches the scenario where consistency breaks down. The key idea of the proof relies on the following sensitivity inequality resulting from Assumption \ref{marginal hypothesis} which holds for every $t\geq 0$ and mapping $w\in W_P$:
\begin{align*}
        \E_{(x,c)\sim P}\left(\left|\hat{c}_\theta (x)^\top  w(\hat{c}_\theta (x)) - \hat{c}_\theta (x)^\top  w(x)\right|\right)\geq tA_1 \norm{\theta}\E_{(x,c)\sim P}\left(\norm{w(\hat{c}_\theta (x)) - w(x)}\right) - A_2\norm{\theta}t^\alpha,
\end{align*}
where $A_1,A_2$ are positive constants. The left-hand side equals $\ell_P^\beta(\theta)$ for a well-chosen $w \in W_P$, while the right-hand side is directly related to $\ell_P(\theta) - \beta$. Optimizing over $t$ in the right-hand side yields the desired inequality. By omitting the term $A_2 | \theta | t^\alpha$ and the dependence of the right-hand side on $t$, it is possible to obtain a stronger bound. Specifically, if we assume the existence of $B_s > 0$ such that for every $\theta \in \R^m$ and every mapping $w \in W_P$, 
\begin{align} \label{no polyhedron inequality}
        \E_{(x,c)\sim P}\left(\left|\hat{c}_\theta (x)^\top  w(\hat{c}_\theta (x)) - \hat{c}_\theta (x)^\top  w(x)\right|\right)\geq B_{s} \norm{\theta}\E_{(x,c)\sim P}\left(\norm{w(\hat{c}_\theta (x)) - w(x)}\right),
\end{align}
the bound $\ell_P(\theta) \leq \beta + O\left(\ell_P^\beta(\theta)\right)$ holds when $\norm{\theta}$ is bounded away from $0$. This bound is stronger since $1-\frac{1}{1+\alpha}<1$ for every $\alpha>0$.
\begin{Remark}
This theorem suggests that to minimize \(\ell_P\), one should seek an approximate minimizer \(\theta\) of \(\ell_{P}^{\beta}\) for a well-chosen $\beta \in \R$ with \(\|\theta\|\) kept away from zero. By applying Theorem \ref{surrogategen}, it follows that it suffices to find some $\theta \in\R^m$ such that \(\ell_{P_n}^{\beta}(\theta)\) is small while ensuring \(\|\theta\|\) remains bounded away from zero.
\end{Remark}

If inequality (\ref{no polyhedron inequality}) holds, we can obtain the bound $\ell_P(\theta) \leq \beta + O\left(\ell_P^\beta(\theta)\right)$ without requiring $W$ to be a polyhedron. 

\subsection{Optimizing our surrogate}
Based on theorems \ref{mainconsistency}, \ref{surrogategen}, and \ref{sensitivity}, it is sufficient to find a \(\theta \in \mathbb{R}^m\) such that \(\ell_{P_n}^{\beta}\) is small and \(\theta\) is kept bounded away from zero to minimize \(\ell_P\).
The empirical surrogate loss $\ell_{P_n}^{\beta}$ can be written as:
\begin{align*}
\forall \theta \in \R^m,\; \ell_{P_n}^\beta (\theta)=g_{P_n}(\theta) - \overline{g}_{P_n}^\beta (\theta),
\end{align*}
where
\begin{eqnarray}
&&g_{P_n}(\theta):=-\mathbb E_{(c,x)\sim {P_n}}\left(\min_{w\in W}\hat{c}_{\theta}(x)^\top w\right)=-\min_{{w}_{P_n}\in W_{P_n}}\mathbb E_{(c,x)\sim {P_n}}\left(\hat{c}_{\theta}(x)^\top w_{P_n}(x)\right) \label{min function equality}\\
&&\overline{g}_{P_n}^{\beta}(\theta):=-\min_{{w}_{P_n}\in \overline{W}_{P_n}^{\beta}}\mathbb E_{(c,x)\sim {P_n}}\left(\hat{c}_{\theta}(x)^\top \overline{w}_{P_n}^\beta(x)\right).
\end{eqnarray}
The equality in (\ref{min function equality}) is satisfied because choosing a policy $w_P$ minimizing $E_{(c,x)\sim {P_n}}\left(\hat{c}_{\theta}(x)^\top \overline{w}_{P_n}^\beta(x)\right)$ is equivalent to choosing $w_P(x)$ minimizing $\hat{c}_\theta(x)^\top w_P(x)$ for every $x$ in $\R^k$.
\begin{paragraph}{\textbf{Landscape properties and smoothing.} } We take a closer look at the structure of our surrogate. We make the following mild assumption.

\begin{Assumption} \label{smoothness}
    The function $\theta \longmapsto \hat{c}_\theta(x)$ is differentiable for every $x\in \R^k$, and there exists $B_L>0$ such that for all $\theta,\theta'\in \R^m$ and $x\in \R^k$,
    \begin{align*}
        \norm{\nabla_\theta\hat{c}_\theta(x) - \nabla_\theta\hat{c}_{\theta'}(x)}\leq B_L \norm{\theta - \theta'}.
    \end{align*}
\end{Assumption}
This assumption is sufficient for $\ell_{P_n}^\beta$ to be written as a difference of two convex functions. Specifically, we have the following result:
\begin{Proposition} \label{DC structure}
    Under assumptions \ref{smoothness} and \ref{boundedness}, $\ell_{P_n}^\beta$ is a DC (difference of two convex functions) function.
\end{Proposition}
The proof of the proposition above (see appendix \ref{DC property proof}) relies on the fact that under assumptions \ref{smoothness} and \ref{boundedness}, $g_{P_n}$ and $\overline{g}_{P_n}^\beta$ are both weakly convex functions, and consequently their difference is a difference of convex functions.

We aim to find $\theta \in \R^m \setminus \{0\}$ that minimizes $\ell_{P_n}^\beta$. A natural approach is to identify a stationary point of $\ell_{P_n}^\beta$ (i.e., a $\theta$ where the subgradient of $\ell_{P_n}^\beta$ contains 0). Despite $\ell_{P_n}^\beta$ having a DC structure, it may still be non-smooth, making the task of finding a stationary point challenging. To address this, we introduce a smoothed version of $\ell_{P_n}^\beta$.

\begin{Definition}{(s-CILO loss)} \label{sCILO def}
    For all $\lambda \in \R^m$ and $\beta \geq \beta_{\mathcal H, P_n}$, let
    \begin{align*}
        M_{P_n}(\lambda):=\min_{\theta\in \R^m}(g_{P_n}(\theta)+\frac{1}{2}\|\lambda-\theta\|^2),\;\overline{M}_{P_n}^{\beta}(\lambda):=\min_{\theta\in \R^m}(\overline{g}_{P_n}^{\beta}(\theta)+\frac{1}{2}\|\theta-\lambda\|^2) 
    \end{align*}
    \vspace{-4mm}
    \begin{align*}
    \theta_{P_n}(\lambda)=\arg\min_{\theta \in \R^m}(g_{P_n}(\theta)+\frac{1}{2}\|\theta-\lambda\|^2), \text{ and }\overline{\theta}_{P_n}^{\beta}(\lambda)=\arg\min_{\theta\in \R^m}(\overline{g}_{P_n}^{\beta}(\theta)+\frac{1}{2}\|\theta-\lambda\|^2).
    \end{align*}
    
We define a smooth surrogate to the CILO loss, which we call the s-CILO loss
\begin{equation}
r_{P_n}^{\beta}(\lambda):=M_{P_n}(\lambda)-\overline{M}_{P_n}^{\beta}(\lambda).
\end{equation}
\end{Definition}
Following \cite{sun2021algorithms} (proposition 1, page 10), the smooth surrogate above has the following property: for every stationary point $\lambda$ of $r_{P_n}^\beta$, $\overline{\theta}_{P_n}^\beta (\lambda)$ and $\theta_{P_n}(\lambda)$ are equal and are stationary points of $\ell_{P_n}^\beta$.
Similar to $\ell_{P_n}^\beta$, s-CILO satisfies the positivity property and the fact that its minimum value is equal to $0$ from $\ell_{P_n}^\beta$. 

We now have a way to find a stationary point of $\ell_{P_n}^\beta$.  
A practical way to ensure that the iterates of an optimization algorithm for $\ell_{P_n}^\beta$ do not converge to $\theta=0$ is to add a constraint $\norm{\theta}\geq z$ where $z>0$.
 Since the function $\theta \longmapsto z-\norm{\theta}$ can be written as a difference of two convex functions, we can find a constrained stationary point to $r_{P_n}^\beta$ (see \cite{pang2017computing} for tractability results for DC constrained minimization). 

The DC structure of our surrogate enables us to find a stationary point of $\ell_{P_n}^\beta$, thus allowing us to approximately minimize the target loss $\ell_P$. Moreover, as we will demonstrate (see Theorem \ref{subdifferentialsurrogate} below), for specific choices of the hypothesis set (e.g., a linear hypothesis set), our surrogate has a favorable landscape, meaning that every stationary point of $\ell_{P_n}^\beta$ is a global minimum. Consequently, first-order optimization algorithms can effectively find a global minimum of our surrogate.

From now on, we assume that our hypothesis set $\mathcal H$ is linear, meaning there exists a mapping $\Phi: \R^k \longrightarrow \R^{m \times d}$ such that for all $\theta \in \R^m$, $\hat{c}_\theta(x) = \Phi(x) \theta$. In this scenario, Assumption \ref{characterize zero} is satisfied if $\Phi(x)$ is full rank almost surely. Additionally, Assumption \ref{smoothness} is inherently satisfied, and Assumption \ref{normboundedness} can be replaced with the following:

\begin{Assumption} \label{feature matrix bounded} For any $x$, the largest singular value of $\Phi(x)$ is bounded above by $B_{\Phi}$. \end{Assumption}

With this special choice of hypothesis set,
our surrogate loss and its smoothed version enjoy good landscape properties. In particular, every stationary point of $\ell_{P_n}^\beta$ is a global minimum. We formulate our main landscape result.
\begin{Theorem} \label{subdifferentialsurrogate}
\;
\begin{enumerate}
\item For all $\theta \in \R^m$, $\overline{g}_{P_n}^\beta$ and $g_{P_n}$ have a non-empty subgradient at $\theta$, and $\ell_{P_n}^{\beta}(\theta)=\theta^\top v(\theta)$ for some $v(\theta)\in \partial \ell_{P_n}^\beta:=\partial\overline{g}_{P_n}^{\beta}(\theta) -\partial g_{P_n}^{\beta}(\theta):=\{u_\beta - u,\; (u_\beta,u)\in \partial\overline{g}_{P_n}^{\beta}(\theta) \times \partial g_{P_n}(\theta)\}$, where $\partial \overline{g}_{P_n}^{\beta}(\theta)$ and $\partial g_{P_n}(\theta)$ are respectively the subgradients of  $\overline{g}_{P_n}^\beta$ and $g_{P_n}$ at $\theta$;
\item For all $\theta \in \R^m$, if $v(\theta)=0$ for some $v(\theta)\in \partial \ell_{P_n}^{\beta}(\theta)$, then $\ell_{P_n}^{\beta}(\theta)=0$, hence $\theta$ is a minimizer.
\item $r_{P_n}^\beta$ is everywhere differentiable. If $\lambda$ is a stationary point of $r_{P_n}^{\beta}$ such that $\hat{c}_{\overline{\theta}_{P_n}^{\beta}(\lambda)}(x)\neq 0$ almost surely, then $\overline{\theta}_{P_n}^{\beta}(\lambda)$  is a global minimum of $\ell_{P_n}^{\beta}$. In particular, under assumptions \ref{boundedness} and \ref{normboundedness}, we have for every $\theta \in \R^m$ and $\varepsilon \geq 0$,
    \begin{align*}
        \norm{\nabla r_{P_n}^\beta}\leq \varepsilon \Longrightarrow \ell_{P_n}^\beta(\theta)\leq 8B_WB_\Phi \varepsilon. 
    \end{align*}
\end{enumerate}
\end{Theorem}
The proof of this theorem is in appendix \ref{subdifferentialsurrogate proof}. 

Assuming a priori that running gradient descent to optimize our smooth surrogate, we do not fall into the case $\hat{c}_\theta(x)=0$ almost surely, our practical procedure to optimize $\ell_P$ is the following. In the algorithm below, $P_n$ is the uniform distribution over the dataset $S_{\text{train}}$, and $\hat P$ is the uniform distribution over the testing dataset $S_\text{test}$.

\begin{algorithm}
\caption{Algorithm to optimize $\ell_P$}
\begin{algorithmic}
\STATE \textbf{Input:} Training and testing datasets $S_{\text{train}}$, $S_{\text{test}}$ drawn from $P$, line search precision level $\eta$, $\overline \beta,\; \underline \beta$ upper and lower bounds for $\beta$.
\STATE $ \Theta\leftarrow \varnothing.$
\FOR{$i = 0$ to $\left\lfloor \frac{\overline \beta - \underline \beta}{\eta} \right\rfloor$}
    \STATE $\beta \leftarrow \underline \beta + i \times \eta$
    \STATE Run gradient descent on $r_{P_n}^\beta$.
    \STATE Set $\lambda$ to be the gradient descent iterate which yields the smallest value of $\ell_{\hat P}(\overline\theta_{P_n}^\beta(\lambda))$.
    \STATE Add $\overline\theta_{P_n}^\beta(\lambda)$ to $\Theta$.
\ENDFOR
\RETURN $\arg\min_{\theta \in \Theta} \ell_{\hat P}(\theta)$.
\end{algorithmic}
\end{algorithm}

\end{paragraph}
Now, the only remaining problem to address is the case where the algorithm gives $\theta=0$.
\begin{paragraph}{\textbf{Avoiding zero solutions.} } 
Even if we successfully minimize \(\ell_{P_n}^\beta\), we could encounter the pathological case where \(\theta = 0\), which would fail to ensure that \(\hat{c}_{\theta}(x) \neq 0\) almost surely. In this situation, the conditions for Theorems \ref{mainconsistency} and \ref{sensitivity} would not be satisfied. To prevent this, we leverage the linear structure of the hypothesis set. We begin by stating a few propositions that highlight key properties of our surrogate when the hypothesis set is linear.
\begin{Definition}
For a set $V$ in $\R^m$ and element $u\in \R^m$, 
denote $-V=\{-v,\; v\in V\} \text{ and } d(u,V)=\min_{v\in V}\norm{u-v}$ the $L^2$ distance between $u$ and $V$.
\end{Definition}
\begin{Proposition}\label{representation}
Let $\lambda \in \R^m$ and $\beta \geq \beta^\star_{\mathcal H,P_n}$, we denote
\begin{align*}
    V_{P_n}:=\left\{\E_{(x,c)\sim P}\left(\Phi^\top (x)w_{P_n}(x)\right),w_{P_n}\in W_{P_n}\right\},\;\overline{V}^\beta_{P_n}:=\left\{\E_{(x,c)\sim P}\left(\Phi^\top (x)w_{P_n}(x)\right),w_{P_n}\in \overline{W}^\beta_{P_n}\right\}.
\end{align*}
$r_{P_n}^\beta$ can be rewritten as 
\begin{align*}
    r_{P_n}^\beta(\lambda)=\displaystyle\min_{v \in \overline{V}_{P_n}^\beta}\frac{1}{2}\norm{\lambda + v}^2 - \displaystyle\min_{v \in V_{P_n}}\frac{1}{2}\norm{\lambda + v}^2=\frac{1}{2}\left(d(\lambda,-\overline{V}_{P_n}^\beta)^2-d(\lambda,-V_{P_n})^2\right)
\end{align*}
\end{Proposition}
The following is a key lemma to prove Proposition \ref{representation}, and also to design the procedure which will enable us to optimize $\ell_{P_n}^{\beta}$.
\begin{Lemma} \label{transferoptimality}
    For every $\lambda \in \R^m$, we have
    \begin{enumerate}
\item $\theta_{P_n}(\lambda)=\lambda+v$ with $v=\arg \displaystyle\min_{v \in V_{P_n}}\frac{1}{2}\norm{\lambda + v}^2$ and $M_{P_n}(\lambda)=\frac{1}{2}\norm{\lambda}^2 - \displaystyle\min_{v \in V_{P_n}}\frac{1}{2}\norm{\lambda + v}^2$;
\item $\overline{\theta}^\beta_{P_n}(\lambda)=\lambda+v_\beta$ with $v_\beta=\arg \displaystyle\min_{v \in \overline{V}_{P_n}^\beta}\frac{1}{2}\norm{\lambda + v}^2$ and $\overline{M}_{P_n}^\beta(\lambda)=\frac{1}{2}\norm{\lambda}^2 - \displaystyle\min_{v \in \overline{V}_{P_n}^\beta}\frac{1}{2}\norm{\lambda + v}^2$.
\end{enumerate}
    \end{Lemma}
The key idea behind the lemma above is to simply swap min and max in the definitions of $M_{P_n}$ and $\overline{M}_{P_n}^\beta$ (see Definition \ref{sCILO def}). The proof of the proposition and the lemma can be found in appendix \ref{representation proof} and \ref{transferoptimality proof}.


We observe that \( r_{P_n}^\beta \) exhibits an interesting structure: when \(\lambda\) minimizers $r_{P_n}^\beta$, its distances from \( -V_{P_n} \) and \( -\overline{V}_{P_n}^\beta \) are equal. However, when running gradient descent on $r_{P_n}^\beta$, we want to avoid the scenario where \( d(\lambda, -\overline{V}_{P_n}^\beta) = d(\lambda, -V_{P_n}) = 0 \). This happens when the gradient descent iterates converge inside the set \( -\overline{V}_{P_n}^\beta \). To address this, we aim to push the gradient descent iterates outside the sets \( -V_{P_n} \) and \( -\overline{V}_{P_n}^\beta \). We achieve this by leveraging the structure of \( r_{P_n}^\beta \), applying the logarithm to the two distances to ensure they do not equal zero.

    \begin{Definition}{(log-CILO loss)} Let $\beta \geq \beta^\star_{\mathcal H,P}$, we define for every $\theta \in \R^m$,
\begin{equation} \label{CILOloss}
        f_{P_n}^{\beta}(\lambda)=\log d(\lambda,\overline V_{P_n}^\beta) - \log d(\lambda, V_{P_n}).
\end{equation}
We call this function the log-CILO loss.

\end{Definition}

The function \( f_{P_n}^\beta \) inherits the properties we previously observed in \( r_{P_n}^\beta \). Specifically, since \( d(\lambda, -\overline{V}_{P_n}^\beta) \geq d(\lambda, -V_{P_n}) \) for every \(\lambda \in \mathbb{R}^m\) (due to the fact that \( \overline{V}_{P_n}^\beta \subset V_{P_n} \)), and because the logarithm is a non-decreasing function, we have \( f_{P_n}^\beta(\lambda) \geq 0 \) for every \(\lambda \in \mathbb{R}^m\). Furthermore, similar to $r_{P_n}^\beta$, if \(\lambda\) is a minimizer of \( f_{P_n}^\beta \), then the distances \( d(\lambda, -\overline{V}_{P_n}^\beta) \) and \( d(\lambda, -V_{P_n}) \) are equal. Additionally, it can be shown that any stationary point of \( f_{P_n}^\beta \) is also a stationary point of \( r_{P_n}^\beta \) (see Theorem \ref{move}).

If gradient descent is used to optimize $ r_{P_n}^\beta $ and a subsequence of the iterates converges outside the set $-V_{P_n}$ (and hence outside $-\overline V_{P_n}^\beta$ as well), we can recover from the limit of this sequence a $\theta \in \R^m$ that is bounded away from zero and serves as a minimizer of $\ell_{P_n}^\beta$. As a result, by Theorems \ref{surrogategen} and \ref{sensitivity}, we successfully obtain a desirable approximate minimizer of $\ell_P$. However, if every limit point of the gradient descent iterates converges within the sets $-V_{P_n}$ and $-V_{P_n}^\beta$, the function $f_{P_n}^\beta$ can be used instead. Gradient descent applied to $f_{P_n}^\beta$ will not yield iterates converging inside the sets $-V_{P_n}$ and $-V_{P_n}^\beta$ due to the presence of $\log$ barriers. The following proposition outlines the possible outcomes when running gradient descent on $f_{P_n}^\beta$.

\begin{Proposition}
Consider the backtracking line search gradient descent method applied to the function \( f_{P_n}^{\beta} \), ensuring sufficient decrease at each iteration (refer to Algorithm 3.1 in \cite{nocedal1999numerical}). Suppose the gradient descent sequence has at least one limit point. There are two possible cases:
\begin{enumerate}
    \item The limit point lies strictly outside the interior of \( -V_{P_n} \) or \( -\overline{V}_{P_n}^\beta \);
    \item The limit point lies on the common boundary of \( -V_{P_n} \) and \( -\overline{V}_{P_n}^\beta \).
\end{enumerate}
\end{Proposition}

\begin{Proofe}
Let \(\{\lambda_i\}\) be the sequence of iterates generated by applying gradient descent to \(f_{P_n}^\beta\). Suppose a subsequence \(\{\lambda_{i_t}\}\) of this sequence converges to a limit point \(\lambda^\star\). If \(\lambda^\star\) is not on the boundary of \( -V_{P_n} \) or \( -\overline{V}_{P_n}^{\beta} \), then by \cite{nocedal1999numerical}, \(\lambda^\star\) is a stationary point of \( f_{P_n}^{\beta} \). We now show that if the limit point \(\lambda^\star\) lies on the boundary of \( -V_{P_n} \), then it must also lie on the boundary of \( -\overline{V}_{P_n}^{\beta} \). Indeed, on the one hand, the objective function \( f_{P_n}^{\beta}(\lambda_{i_t}) \) is non-increasing, so \( f_{P_n}^{\beta}(\lambda_{i_t}) \leq f_{P_n}^{\beta}(\lambda_0) \) for all \( t \). On the other hand, \( -\log(d(\lambda^{i_t}, -V_{P_n})) \to \infty \), which implies that \( -\log(d(\lambda^{i_t}, -\overline{V}_{P_n}^{\beta})) \to \infty \). Therefore, \( d(\lambda^\star, -\overline{V}_{P_n}^{\beta}) \to 0 \), yielding the desired result.
\end{Proofe}

If the gradient sequence \(\{\lambda_{i_t}\}\) converges to the common boundary of \( -V_{P_n} \) and \( -\overline{V}_{P_n}^\beta \), the corresponding sequence \(\{\overline{\theta}_{P_n}^\beta(\lambda_{i_t})\}\) approaches zero, which is not a desirable solution. Fortunately, this issue can be addressed. Since the sequence converges to a limit point on the boundary of \( -V_{P_n} \) and \( -\overline{V}_{P_n}^{\beta} \) which is also a stationary point of $r_{P_n}^\beta$, there must exist some \(\lambda_{i_t}\) such that is an \(\varepsilon\)-stationary point of \( r_{P_n}^{\beta} \). We then move \(\lambda_{i_t}\) orthogonally away from the boundary of \( -V_{P_n} \) (see Figure \ref{fig:enter-label}). This yields a new approximate stationary point that is sufficiently distant from both \( -V_{P_n} \) and \( -V_{P_n}^\beta \). As a result, the corresponding \(\theta\) is an approximate global minimizer of \(\ell_{P_n}^{\beta}\) that is bounded away from zero.

The following theorem provides an overview of this procedure, as well as a relationship between the optimality for $f_{P_n}^\beta$ and $r_{P_n}^\beta$.

\begin{Theorem}\label{move}
We define $F_{V_{P_n}}(\lambda):=\arg\min_{\lambda'\in -V_{P_n}}\norm{\lambda - \lambda'}$ and $F_{\overline V_{P_n}^\beta}(\lambda):=\arg\min_{\lambda'\in -\overline V_{P_n}^\beta}\norm{\lambda - \lambda'}$ to be respectively the projection of $\lambda$ on $-V_{P_n}$ and on $-\overline V_{P_n}^\beta$. Suppose that $\lambda$ is an $\varepsilon$-solution of $f_{P_n}^\beta$, i.e. $\norm{\nabla f_{P_n}^\beta (\lambda)}\leq \varepsilon$ for some $\varepsilon >0$, and $\lambda\notin V_{P_n}^\beta$. We have $\norm{\nabla r_{P_n}^\beta}\leq 27 B_W^2B_\Phi^2 \varepsilon$. Furthermore, if $d(\lambda,-V_{P_n})\leq B_WB_\Phi$,  letting  $\tilde{\lambda}=\lambda + r \frac{\lambda-F_{V_{P_n}}(\lambda)}{\|\lambda-F_{V_{P_n}}(\lambda)\|}$ with $r\in [3B_WB_\Phi,8B_WB_\Phi]$ and $\varepsilon'=27 B_W^2B_\Phi^2 \varepsilon$, we have

\begin{itemize}
\item The value of $r_{P_n}^\beta(\tilde{\lambda})$ is dominated by $\varepsilon$. In particular, we have $r_{P_n}^\beta(\Tilde{\lambda})\leq 11B_WB_\Phi \varepsilon'$ and $\norm{\nabla r_{P_n}^\beta (\Tilde{\lambda})}\leq \sqrt{11B_WB_\Phi \varepsilon'}$.
\item The norm of the resulting candidate solution for $\ell_{P_n}^\beta$ is bounded from above and below. In particular, $B_WB_\Phi\le \|\bar{\theta}_{P_n}^\beta(\tilde{\lambda})\|\le 11B_WB_\Phi$.
\end{itemize}
\end{Theorem}

Figure \ref{fig:enter-label} summarizes the construction of $\Tilde{\lambda}$.

\begin{center}
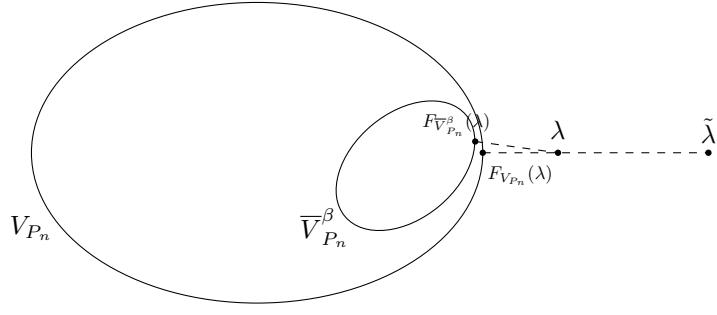
\begin{figure}
    \centering
    \begin{tikzpicture}

\begin{scope}[rotate=40, scale=0.7, shift={(0.5,0)}]
    \draw (1.5,-2) ellipse (1.5cm and 1cm);
    
\end{scope}
\node at (-3, -1) {$V_{P_n}$};
\draw (0,0) ellipse (3cm and 2cm);
\node at (0.9, -1) {$\overline{V}^{\beta}_{P_n}$};

\draw[dashed, -] (3,0) -- (6,0);
\draw[dashed, -] (2.9,0.15) -- (4,0);
\node at (4, 0.3) {$\lambda$};
\node at (6, 0.3) {$\Tilde{\lambda}$};
\node[scale=0.7] at (3.5, -0.3) {$F_{V_{P_n}}(\lambda)$};
\node[scale=0.7] at (2.65, 0.4) { $F_{\smash{\overline V_{P_n}^\beta}}(\lambda)$};
\fill (4,0) circle (1.2pt);
\fill (6,0) circle (1.2pt);
\fill (3,0) circle (1.2pt);
\fill (2.9,0.15) circle (1.2pt);


\end{tikzpicture}
    \caption{Construction of $\Tilde{\lambda}$: when $\lambda$ is close to the sets $\overline{V}_{P_n}^\beta$ and $V_{P_n}$, we construct a new approximately stationary point by moving away from $V_{P_n}$ orthogonally to its boundary.}
    \label{fig:enter-label}
\end{figure}
\end{center}

Hence, the theorem above provides us with an algorithm to optimize $\ell_{P_n}^\beta$:
\begin{enumerate}
    \item Run gradient descent on $r_{P_n}^\beta$.
    \item If we obtain a solution reasonably far from $V_{P_n}$, we simply return that iterate as an output. Else, we run gradient descent on $f_{P_n}^\beta$.
    \item If the gradient descent iterates do not converge to a common boundary to $-V_{P_n}$ and $-\overline V_{P_n}^\beta$, we return the final iterate as an output. Else, we use the procedure in Theorem \ref{move} to obtain a near stationary point which is far from the boundary, and return it as an output.
\end{enumerate}
\end{paragraph}

We can now combine all our results to provide optimality guarantees for \(\ell_P\). The procedure described above enables us to find a limit point $\lambda$ that is an \(\varepsilon\)-stationary solution for \(r_{P_n}^\beta\). As a result, by Theorem \ref{subdifferentialsurrogate} and Theorem \ref{move}, we obtain a near-stationary point for \(\ell_{P_n}^\beta\) that is bounded away from zero. Furthermore, Theorem \ref{surrogategen} allows us to conclude that we have achieved a good solution for \(\ell_P^\beta\). Theorem \ref{sensitivity} confirms that we have indeed found a solution \(\theta \in \mathbb{R}^m\) such that \(\ell_P(\theta) \leq \beta\), up to a small approximation error. Finally, choosing a suitable $\beta$ using line search provides us with a solution $\theta^\star$ such that $\ell_P(\theta^\star)$ is small.
\section{Computational experiments} \label{sec4}
\subsection{Experimental setting}
To create a framework where we can enforce misspecification and progressively analyze how the performance of different methods varies with the level of misspecification $\gamma_{\text{miss}}(\mathcal{H})$, we need to define a model that allows for easy adjustment of the misspecification level while still yielding meaningful results. We define a model which can easily be slightly changed to be misspecified. In order to do so, we define for $x \in  \R^k$, $c(x)$ to be a linear combination of polynomial functions of $x$. This could be written in a compact way as $c(x)=M\phi(x)$, where $M$ is a real valued matrix, and $\phi(x)$ is a vector whose coordinates are polynomial in $x$. Furthermore, we define $\hat{c}_\theta (x)$ for $x \in \R^k$ and $\theta \in \R^m$ as a linear combination of the coordinates of $\phi(x)$. To increase the level of misspecification, we can omit some coordinates of $\phi(x)$ in the linear combination. This can also be written compactly as $\hat{c}_\theta(x) = M(\theta) \Tilde{\phi}(x)$, where $M(\theta)$ is a matrix representation of $\theta$ and $\Tilde{\phi}(x)$ is a truncated version of $\phi(x)$, depending on the desired level of misspecification. We set $W$ to be a polyhedron and written as $W=\{w \in \R^d,\; Aw=b,\; 10\geq w\geq 0\}$ where $A\in \R^{j\times d }$ ($j\leq d$) and $b\in \R^j$. $W$ is closed and bounded, so Assumption \ref{boundedness} holds. In every experiment, we sample $x\sim \mathcal N(0,I)$ while all of its coordinates are conditioned to be between $0$ and $10$. and the coefficients of $A$ from a standard normal Gaussian distribution, and $b$ to be equal to $A\abs{w}$, where $w$ has standard normal random coefficients. We set $(d,k)=(20,5)$. We write 
\begin{align*}
    \forall x=(x_1,\dots,x_5) \in \R^5,\; \phi(x)=\begin{pmatrix}
        x_1 & \dots & x_5 & x_1x_2 & \dots & x_1x_2x_3x_4x_5
    \end{pmatrix}^\top =\left(\prod_{i=1}^{5}x_i^{y_i}\right)_{y\in \{0,1\}^5,y\neq 0}.
\end{align*}
Furthermore, for any matrix $M$, we denote $L_r(M)$ the $r-$th row of $M$. We set the ground truth model to be $c(x)=M(\theta^\star)\phi(x)$ where $\theta^\star=\begin{pmatrix}
    L_1(M(\theta^\star)) & \dots & L_5(M(\theta^\star))
\end{pmatrix}^\top \in \R^{155}$. We also write for $\theta \in \R^{m}$, $\hat{c}_\theta (x) = M(\theta)\Tilde{\phi}(x)$, where $m=5\times (31-s)$ and  $\theta=\begin{pmatrix}
    L_1(M(\theta)) & \dots & L_5(M(\theta)) 
\end{pmatrix}^\top \in \R^{5(31-s)}$. Here, $s$ is the number of coordinates (or features) we will be removing from $\phi(x)$ to obtain $\Tilde{\phi}(x)$, and hence causing the model to be misspecified. This can be written coherently with our theoretical setting. Indeed, we can write for all $x\in \R^k$, $\hat{c}_\theta (x) = \Tilde{\Phi}(x) \theta$ and $c(x)=\Phi(x)\theta^\star$, $\theta^\star \in \R^m$, where 
\begin{align*}
    \Tilde{\Phi}(x) = \begin{pmatrix}
        \Tilde{\phi}(x)^\top  & & 0\\
        & \ddots & \\
        0 & & \Tilde{\phi}(x)^\top 
    \end{pmatrix} \in \R^{20 \times (5(31-s))} \text{ and }\Phi(x) = \begin{pmatrix}
        \phi(x)^\top  & & 0\\
        & \ddots & \\
        0 & & \phi(x)^\top 
    \end{pmatrix} \in \R^{20 \times 155}.
\end{align*}
We denote $\mathcal H_s=\{\hat{c}_\theta : x \longmapsto \Tilde \Phi (x) \theta,\; \theta \in \R^m\}$ the resulting hypothesis set when the number of removed coordinates from $\phi$ is equal to $s$. We have, for any $s>0,$ $\gamma_{\text{miss}}(\mathcal H_s)>0$. Also, since for any $s,s'\geq 0$ such that $s\geq s'$, the inclusion $\mathcal H_{s'}\subset \mathcal H_s$ is satisfied, we have $\gamma_{\text{miss}}(\mathcal H_s)\geq \gamma_{\text{miss}}(\mathcal H_{s'})$. In other words, $\gamma_{\text{miss}}(\mathcal H_s)$ is increasing with $s$. Furthermore, since the coordinates of the matrices \(\Tilde{\Phi}(x)\) and \(\Phi(x)\) are bounded, Assumptions \ref{boundedness2} and \ref{normboundedness} are satisfied. Additionally, because the function \(\hat{c}_\theta\) is nonzero and analytic for any \(\theta \neq 0\), and \(\hat{c}_\theta(x)^\top v \neq 0\) for any nonzero \(v\) in the polyhedron \(W\), Proposition \ref{continuousdist} ensures that Assumption \ref{regularity} holds. Moreover, Assumption \ref{smoothness} is satisfied because the hypothesis set \(\mathcal H\) is linear.

\subsection{Experimental results}
\begin{center}
\begin{figure}[h]
    \centering
        \includegraphics[scale=0.8]{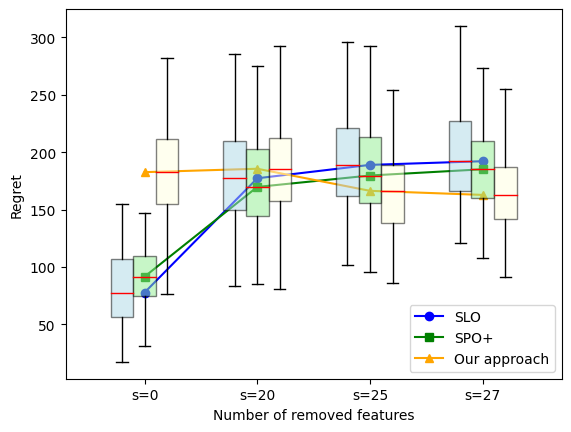}
        
    \caption{SPO+ v.s. SLO v.s. our approach comparison for different levels of misspecification}
    \label{spo_vs_slo_vs_CILO}
\end{figure}
\end{center}

We revisit the comparison performed by \cite{hu2022fast}. In their work, SLO is compared to SPO+, showing that when the hypothesis set is well-specified, SLO achieves a smaller value of $\ell_P$ than SPO+. Conversely, when the model is misspecified, SPO+ performs better. We extend this comparison by including our method as well. The exact losses minimized for each method are as follows:
\begin{align*}
    \ell_{\text{SLO}}(\theta)=\mathbb{E}_{(x,c)\sim P_n}\left(\norm{\hat{c}_\theta (x) - c(x)}^2\right),\; \ell_{P_n}^+(\theta):=\E_{(x,c)\sim P_n}\left(\ell_ {\text{SPO+}}(\hat{c}_{\theta}(x),c)\right) \\
    \ell_{P_n}^{\beta}(\theta)= \min_{\overline{w}_{P_n}^{\beta}\in \overline{W}_{P_n}^{\beta}}\E_{(x,c)\sim P_n}\left(\hat{c}_\theta^\top \overline{w}_{P_n}^{\beta}(x)\right) - \min_{w_{P_n}\in W_{P_n}}\E_{(x,c)\sim P_n}\left(\hat{c}_\theta^\top w_{P_n}(x)\right).
\end{align*}
To optimize $\ell_{P_n}^\beta$, we run gradient descent on its surrogate loss $r_{P_n}^\beta$ (we did not need to optimize the surrogate $f_{P_n}^\beta$ because the iterates did not converge to $0$). We choose $\beta$ by line search. We use $\beta_{\text{min},P}=\E_{(x,c)\sim P_n}\left(c^\top w(c)\right)$ as a lower bound to $\beta$, and $\beta_{\text{SPO+}}=\E_{(x,c) \sim P_n}\left(c^\top w\left(\hat{c}_{\theta^\star_{\text{SPO+}}}(x)\right)\right)$ where $\theta^\star_{\text{SPO+}}$ is the solution obtained by optimizing $\ell_{P_n}^+$. We test 96 evenly spaced values of $\beta$ in the interval $[\beta_{\text{min},P},\beta_{\text{SPO+}}]$, and pick $\beta$ yielding the solution with the best decision performance. $\beta^\star_{\mathcal H,P}$ is likely to fall into (or near) the interval $[\beta_{\text{min},P},\beta_{\text{SPO+}}]$, and enjoys good optimality guarantees for $\ell_P$, which ensures that the solution we obtain yields a small value for $\ell_P$. For every value of $s\in \{0,20,25,27\}$, we run 96 experiments to get the distribution of the testing loss $\ell_P$ for every method, and get results in the box plots in Figure \ref{spo_vs_slo_vs_CILO}. For every experiment, we sample $A,b$, the ground truth parameter $\theta^*$, and a training dataset containing 20 samples drawn from distribution $P$, and run gradient descent on the three loss functions we are optimizing. All the experiments were run on the MIT Lincoln Lab supercloud (\cite{reuther2018interactive}).
In Figure \ref{spo_vs_slo_vs_CILO}, we plot the regret yielded by SLO, SPO+ and our approach, i.e. the value of 
$$\ell_P'(\theta)=\ell_P(\theta) - \beta_{\text{min,P}}\geq 0.$$
We can see that the more the model is misspecified, our method performs better compared to SPO+ and SLO. This is because whereas it is unclear whether SPO+ and SLO are capable of minimizing $\ell_P$, our approach enjoys good theoretical guarantees. However, the more well-specified the model is, the better is SLO, which is coherent with the conclusions made in \cite{hu2022fast} and \cite{elmachtoub2023estimate}. We also notice that removing features (as long as a reasonable number of features remain) improves the performance of our approach. One possible reason for this is that our approach overfits more easily than other methods when the hypothesis set is complex (and consequently far from the misspecified case). Further investigation is yet to be done to explain this phenomenon. 
\section{Conclusion} \label{sec5}
In this paper, we provide a novel approach to address model misspecification in contextual optimization. State-of-the-art methods are mainly capable of optimizing the target loss $\ell_P$ representing the decision performance resulting from predictions of the random costs only when the chosen hypothesis set is well-specified. In the misspecified setting, it has been unclear whether it is possible to minimize $\ell_P$ and provide theoretical guarantees. Our surrogate loss function successfully optimizes the target loss and retrieves a good element from the chosen hypothesis set in terms of decision performance in a reasonable time. Even though our surrogate is non-convex and non-smooth, we exploited its structure: it is a difference of two convex functions, which could be smoothed in order to be optimized. We show experimentally and theoretically that our approach outperforms state-of-the-art methods when the hypothesis set is misspecified. To our knowledge, our approach is the first to provably optimize the target loss when the hypothesis set is misspecified, without requiring regularity of $\ell_P$. Although we have proven theoretical guarantees when the chosen hypothesis set is linear, we strongly believe that these results can be generalized to a wider class of predictors. Furthermore, we have run experiments using synthetic data, and have yet to do more experiments using real world datasets to compare our approach to other state-of-the-art methods.

\bibliography{main}
\bibliographystyle{apalike}

\newpage
\section*{Appendix}
\subsection{Proof of Theorem \ref{surrogategen}}
\begin{Proofe}
We proceed in two steps.

\textbf{Step 1.} We first prove the Lipschitz continuity of the generalization, i.e. that there exists a constant $A\geq 0$ such that for every $\theta,\theta'\in \R^m$, \begin{align}\left|\ell_{P_n}^{\beta}(\theta) - \ell_{P}^{\beta}(\theta)-\left(\ell_{P_n}^{\beta}(\theta') - \ell_{P}^{\beta}(\theta')\right)\right|\leq A\norm{\theta - \theta'}.\label{gridineq}\end{align}
  For every $\theta\in \R^m$, using Lagrangian duality,
    \begin{align}
        \ell_{P_n}^{\beta}(\theta)&=\min_{\overline{w}_{P_n}^\beta \in \overline{W}_{P_n}^\beta}\E_{(x,c)\sim P_n}\left(\hat{c}_{\theta}(x)^\top \overline{w}_{P_n}^\beta(x)\right) - \min_{w_{P_n} \in W_{P_n}}\E_{(x,c)\sim P_n}\left(\hat{c}_{\theta}(x)^\top w_{P_n}(x)  \right)\\
        &=\min_{\overline{w}_{P_n}^\beta \in W_{P_n}}\max_{y\geq 0}\E_{(x,c)\sim P_n}\left(\hat{c}_{\theta}(x)^\top \overline{w}_{P_n}^\beta(x)\right) + y \left(\E_{(x,c)\sim P_n}\left(c(x)^\top w_{P_n}^{\beta}(x)\right)- \beta\right) \\&
        - \min_{w_{P_n} \in W_{P_n}}\E_{(x,c)\sim P_n}\left(\hat{c}_{\theta}(x)^\top w_{P_n}(x)  \right)\\
        &=\max_{y\geq 0}\min_{\overline{w}_{P_n}^\beta \in W_{P_n}}\E_{(x,c)\sim P_n}\left(\left(\hat{c}_{\theta}(x)+yc(x)\right)^\top \overline{w}_{P_n}^\beta(x)\right) -y\beta - \min_{w_{P_n} \in W_{P_n}}\E_{(x,c)\sim P_n}\left(\hat{c}_{\theta}(x)^\top w_{P_n}(x)  \right).  \label{switch minmax}
    \end{align}
    In (\ref{switch minmax}), we have switched the min and the max because of strong duality, since the objective function we are optimizing is linear. For a given $\theta \in \R^m$, we denote $y^\star(\theta)$ to be the optimal dual variable corresponding to $\theta$ in the minimization problem above. Let $\theta,\theta'\in \R^m$ and $D\geq \max(\norm{\theta},\norm{\theta'})$. We have 
    \begin{align*}
        \left|\ell_{P_n}^{\beta}(\theta) - \ell_{P}^{\beta}(\theta)-\left(\ell_{P_n}^{\beta}(\theta') - \ell_{P}^{\beta}(\theta')\right)\right| \leq \left|\ell_{P_n}^{\beta}(\theta)-\ell_{P_n}^{\beta}(\theta')\right| + \left|\ell_{P}^{\beta}(\theta)-\ell_{P}^{\beta}(\theta')\right|.
    \end{align*}

    We now bound the two terms above on the right side of the inequality. We have 
    \begin{align*}
        \left|\ell_{P_n}^\beta(\theta)-\ell_{P_n}^\beta(\theta')\right|&=\left|\max_{y\geq 0}\min_{\overline{w}_{P_n}^\beta \in W_{P_n}}\E_{(x,c)\sim P_n}\left(\left(\hat{c}_{\theta}(x)+yc(x)\right)^\top \overline{w}_{P_n}^\beta(x)\right) -y\beta - \min_{w_{P_n} \in W_{P_n}}\E_{(x,c)\sim P_n}\left(\hat{c}_{\theta}(x)^\top w_{P_n}(x)  \right)\right.\\
        & -\left.\max_{y\geq 0}\min_{\overline{w}_{P_n}^\beta \in W_{P_n}}\E_{(x,c)\sim P_n}\left(\left(\hat{c}_{\theta'}(x)+yc(x)\right)^\top \overline{w}_{P_n}^\beta(x)\right) -y\beta - \min_{w_{P_n} \in W_{P_n}}\E_{(x,c)\sim P_n}\left(\hat{c}_{\theta'}(x)^\top w_{P_n}(x)  \right)\right|\\
        &=\left|\min_{\overline{w}_{P_n}^\beta \in W_{P_n}}\E_{(x,c)\sim P_n}\left(\left(\hat{c}_{\theta}(x)+y^\star(\theta)c(x)\right)^\top \overline{w}_{P_n}^\beta(x)\right) -y^\star(\theta)\beta \right. \\& \left. - \min_{w_{P_n} \in W_{P_n}}\E_{(x,c)\sim P_n}\left(\hat{c}_{\theta}(x)^\top w_{P_n}(x)  \right)\right.\\
        & -\left.\min_{\overline{w}_{P_n}^\beta \in W_{P_n}}\E_{(x,c)\sim P_n}\left(\left(\hat{c}_{\theta'}(x)+y^\star(\theta')c(x)\right)^\top \overline{w}_{P_n}^\beta(x)\right) \right. \\&\left.-y^\star(\theta')\beta - \min_{w_{P_n} \in W_{P_n}}\E_{(x,c)\sim P_n}\left(\hat{c}_{\theta'}(x)^\top w_{P_n}(x)  \right)\right|.
    \end{align*}
    In the left hand side of the inequality above, we can see that if $\ell_{P_n}^\beta(\theta)\geq\ell_{P_n}^\beta(\theta')$, then replacing $y^\star(\theta')$ by $y^\star(\theta)$ gives a larger term than the one above, and when $\ell_{P_n}^\beta(\theta)\leq \ell_{P_n}^\beta(\theta')$, replacing $y^\star(\theta)$ by $y^\star(\theta')$ gives a larger term than the one above as well. Let $\Tilde{y}\in \{y^\star(\theta),y^\star(\theta')\}$ such that replacing $y^\star(\theta)$ and $y^\star(\theta')$ by $\Tilde{y}$ makes the expression above larger. We have, by replacing $y^\star(\theta)$ and $y^\star(\theta')$ by $\Tilde{y}$ in the inequality above, we get
    \begin{align*}
         \left|\ell_{P_n}^\beta(\theta)-\ell_{P_n}^\beta(\theta')\right|&\leq \left|\min_{\overline{w}_{P_n}^\beta \in W_{P_n}}\E_{(x,c)\sim P_n}\left(\left(\hat{c}_{\theta}(x)+\Tilde{y}c(x)\right)^\top \overline{w}_{P_n}^\beta(x)\right) - \min_{w_{P_n} \in W_{P_n}}\E_{(x,c)\sim P_n}\left(\hat{c}_{\theta}(x)^\top w_{P_n}(x)  \right)\right.\\
        & -\left.\min_{\overline{w}_{P_n}^\beta \in W_{P_n}}\E_{(x,c)\sim P_n}\left(\left(\hat{c}_{\theta'}(x)+\Tilde{y}c(x)\right)^\top \overline{w}_{P_n}^\beta(x)\right) + \min_{w_{P_n} \in W_{P_n}}\E_{(x,c)\sim P_n}\left(\hat{c}_{\theta'}(x)^\top w_{P_n}(x)  \right)\right|\\
        &= \left|\E_{(x,c)\sim P_n}\left(\min_{w \in W}\left(\hat{c}_{\theta}(x)+\Tilde{y}c(x)\right)^\top w\right) - \E_{(x,c)\sim P_n}\left(\min_{w \in W}\hat{c}_{\theta}(x)^\top w  \right)\right.\\
        & -\left.\E_{(x,c)\sim P_n}\left(\min_{w \in W}\left(\hat{c}_{\theta'}(x)+\Tilde{y}c(x)\right)^\top w\right) + \E_{(x,c)\sim P_n}\left(\min_{w \in W}\hat{c}_{\theta'}(x)^\top w  \right)\right|
        \\&=\left| f_{\Tilde{y},P_n}(\theta)-f_{\Tilde{y},P_n}(\theta') \right|,
    \end{align*}
    where for every $\theta \in \R^m$, 
    $$f_{\Tilde{y},P_n}(\theta)=\E_{(x,c)\sim P_n}\left(\min_{w \in W}\smash{\overbrace{\left(\hat{c}_{\theta}(x)+\Tilde{y}c(x)\right)^\top w}^{h_1(\theta,w,x)}}\right) - \E_{(x,c)\sim P_n}\left(\min_{w \in W}\smash{\overbrace{\hat{c}_{\theta}(x)^\top w}^{h_2(\theta,w,x)}}  \right) $$
    
    $W$ is a convex bounded set, and for all $x\in \R^k$ and the two functions $w\in \R^d$, $\theta \longmapsto h_1(\theta,w,x)$, $\theta \longmapsto h_2(\theta,w,x)$ are both differentiable with respect to $\theta$. Moreover, $\frac{\partial h_1}{\partial \theta}\longmapsto \nabla \hat{c}_\theta (x)^\top w$ and $\frac{\partial h_2}{\partial \theta}\longmapsto \nabla \hat{c}_\theta (x)^\top w$ (where $\nabla\hat{c}_\theta(x)$ is the jacobian of $\theta \longmapsto \hat{c}_\theta(x)$ at $\theta$) are both continuous with respect to $w$ for all $\theta\in \R^m$ and $x\in \R^k$. Hence, using Danskin's theorem, we can say that for all $x\in \R^k$, $\overline{h}_1(\theta,x)=\min_{w\in W}h_1(\theta,w,x)$ and  $\overline{h}_2(\theta,x):=\min_{w\in W}h_2(\theta,w,x)$ are both subdifferentiable, and that for all $\theta\in \R^m$, $\partial \overline{h}_1(\theta,x)=\text{conv}\{\nabla\hat{c}_\theta(x)^\top w,\;w\in \arg\min_{w\in W}h_1(\theta,w,x) \}$, $\partial \overline{h}_2(\theta,x)=\text{conv}\{\nabla\hat{c}_\theta(x)^\top w,\;w\in \arg\min_{w\in W}h_2(\theta,w,x) \}$. These two sets are bounded by $B_WB_\Phi$ because of Assumption \ref{normboundedness}. This means that the two functions in the equality above are both $B_WB_{\Phi}$ lipschitz, which makes $f_{\Tilde{y},P_n}$ $2B_WB_\Phi$ lipschitz.

    We can also easily prove that the inequality above is also true when we replace $P_n$ by $P$. Hence, we can deduce that 
    \begin{align*}
        \left|\ell_{P_n}^{\beta}(\theta) - \ell_{P}^{\beta}(\theta)-\left(\ell_{P_n}^{\beta}(\theta') - \ell_{P}^{\beta}(\theta')\right)\right| &\leq \left|\ell_{P_n}^{\beta}(\theta)-\ell_{P_n}^{\beta}(\theta')\right| + \left|\ell_{P}^{\beta}(\theta)-\ell_{P}^{\beta}(\theta')\right|
        \\&\leq \abs{f_{\Tilde{y},P_n}(\theta) - f_{\Tilde{y},P_n}(\theta')}+\abs{f_{\Tilde{y},P}(\theta) - f_{\Tilde{y},P}(\theta')}
        \\&\leq 4B_WB_\Phi\norm{\theta - \theta'},
    \end{align*}
    which yields the desired result of step 1.

    \textbf{Step 2.} We will now prove the desired generalization bound by bounding the generalization gap locally in balls covering the set which represents the possible values of $\theta^\star$. Let $\mathcal B$ be a set of balls of radius $\gamma>0$ (for the norm $\norm{.}$) such that $\bigcup_{B\in \mathcal B}^{}B=B_{\norm{.}}(0,D)$. According to \cite{wainwright2019high}, it is possible to choose $\mathcal B$ such that $\log |\mathcal B|\leq d\log \left(1+\frac{2D}{\gamma}\right)$. For every $B\in \mathcal B$, let $\theta_{B}$ be an element of $B$. For a given $\varepsilon>0$ we would like to bound the probability $\mathbb P\left(\abs{\ell_{P_n}^{\beta}(\theta) -\ell_{P}^{\beta}(\theta) }>\varepsilon\right)$ using Hoeffding's inequality. In order to do that, we first notice that using Lagrangian duality, we get
    \begin{align*}
        \ell_{P}^\beta(\theta)-\ell_{P_n}^\beta(\theta)&=\max_{y\geq 0}\min_{\overline{w}_{P}^\beta \in W_{P}}\E_{(x,c)\sim P}\left(\left(\hat{c}_{\theta}(x)+yc(x)\right)^\top \overline{w}_{P}^\beta(x)\right) -y\beta - \min_{w_{P} \in W_{P}}\E_{(x,c)\sim P}\left(\hat{c}_{\theta}(x)^\top w_{P}(x)  \right)\\
        &-\max_{y\geq 0}\min_{\overline{w}_{P_n}^\beta \in W_{P_n}}\E_{(x,c)\sim P_n}\left(\left(\hat{c}_{\theta}(x)+yc(x)\right)^\top \overline{w}_{P_n}^\beta(x)\right) -y\beta - \min_{w_{P_n} \in W_{P_n}}\E_{(x,c)\sim P_n}\left(\hat{c}_{\theta}(x)^\top w_{P_n}(x)  \right)\\
        & =\min_{\overline{w}_{P}^\beta \in W_{P}}\E_{(x,c)\sim P}\left(\left(\hat{c}_{\theta}(x)+y_{P}^\star(\theta)c(x)\right)^\top \overline{w}_{P}^\beta(x)\right) -y_{P}^\star(\theta)\beta - \min_{w_{P} \in W_{P}}\E_{(x,c)\sim P}\left(\hat{c}_{\theta}(x)^\top w_{P}(x)  \right)\\ &
        -\left(\min_{\overline{w}_{P_n}^\beta \in W_{P_n}}\E_{(x,c)\sim P_n}\left(\left(\hat{c}_{\theta}(x)+y_{P_n}^\star(\theta)c(x)\right)^\top \overline{w}_{P_n}^\beta(x)\right) - y_{P_n}^\star(\theta)\beta - \min_{w_{P_n} \in W_{P_n}}\E_{(x,c)\sim P_n}\left(\hat{c}_{\theta}(x)^\top w_{P_n}(x)  \right)\right),
    \end{align*}
    where $y^\star_{P}(\theta)$ and $y^{*}_{P_n}(\theta)$ are respectively optimal dual variables in the bottom and top line above. As we have seen before, replacing $y^{*}_{P_n}(\theta)$ by $y^\star_{P}(\theta)$ makes the above term larger. Hence, we have 
    \begin{align*}
          \ell_{P}^\beta(\theta)-\ell_{P_n}^\beta(\theta) &\leq  \min_{\overline{w}_{P}^\beta \in W_{P}}\E_{(x,c)\sim P}\left(\left(\hat{c}_{\theta}(x)+y^\star_{P}(\theta)c(x)\right)^\top \overline{w}_{P}^\beta(x)\right) - \min_{w_{P} \in W_{P}}\E_{(x,c)\sim P}\left(\hat{c}_{\theta}(x)^\top w_{P}(x)  \right) \\& -\left(\min_{\overline{w}_{P_n}^\beta \in W_{P_n}}\E_{(x,c)\sim P_n}\left(\left(\hat{c}_{\theta}(x)+y^\star_{P}(\theta)c(x)\right)^\top \overline{w}_{P_n}^\beta(x)\right) - \min_{w_{P_n} \in W_{P_n}}\E_{(x,c)\sim P_n}\left(\hat{c}_{\theta}(x)^\top w_{P_n}(x)  \right)\right)\\
        &= \underbrace{\E_{(x,c)\sim P}\left(\min_{w\in W}\left(\hat{c}_{\theta}(x)+y^\star_{P}(\theta)c(x)\right)^\top w- \min_{w \in W}\hat{c}_{\theta}(x)^\top w  \right)}_{\Tilde{\ell}_P^{\beta}(\theta)}\\
        &-\underbrace{\E_{(x,c)\sim P_n}\left(\min_{w \in W}\left(\hat{c}_{\theta}(x)+y^\star_{P}(\theta)c(x)\right)^\top w- \min_{w \in W}\hat{c}_{\theta}(x)^\top w  \right)}_{\smash{\Tilde{\ell}^{\beta}_{P_n}(\theta)}}
    \end{align*}
    Hence, we have $\mathbb P\left(\ell_{P}^{\beta}(\theta)-\ell_{P_n}^{\beta}(\theta) >\varepsilon\right) \leq \mathbb P\left( \Tilde{\ell}_{P}^{\beta}(\theta) -\Tilde{\ell}_{P_n}^{\beta}(\theta)>\varepsilon\right)$. This enables us to apply Hoeffding's inequality to the term on the right given that we bound the random variable $$m(x,\theta):=\min_{w \in W}\left(\hat{c}_{\theta}(x)+y^\star_{P}(\theta)c(x)\right)^\top w- \min_{w \in W}\hat{c}_{\theta}(x)^\top w. $$
    It is easy to see that using Cauchy-Schwartz inequality, we have
    \begin{align*}
        \abs{m(x,\theta)}&\leq \left(\norm{\hat{c}_\theta(x)}+\abs{y^\star_{P}(\theta)}\norm{c(x)}\right)B_W + \norm{\hat{c}_\theta(x)}B_W
    \end{align*}
    In order to upper bound the right-hand side, we need to upper bound $y^\star_{P}(\theta)$. In order to do this, we prove that the Lagrange multiplier $y^\star_{P}(\theta)$ is bounded for any $\theta \in \R^m$. Indeed, we have 
    \begin{align}
        \min_{w_{P}\in W_{P}} \E_{(x,c)\sim P}(\hat{c}_\theta (x)^\top  w_{P}(x))&\leq \min_{w_{P}\in W_{P}} \max_{y\geq 0}\E_{(x,c)\sim P}(\hat{c}_\theta (x)^\top  w_{P}(x)) + y (\E_{(x,c)\sim P}(c(x)^\top w_{P}(x))-\beta) \label{explain eval}
        \\
        &=\min_{w_{P}\in W_{P}}\E_{(x,c)\sim P}(\hat{c}_\theta (x)^\top  w_{P}(x)) + y^\star_{P}(\theta) (\E_{(x,c)\sim P}(c(x)^\top w_{P}(x))-\beta)\\&\leq \E_{(x,c)\sim P}(\hat{c}_\theta (x)^\top  w(c(x))) + y^\star_{P}(\theta) (\E_{(x,c)\sim P}(c(x)^\top w(c(x)))-\beta).
    \end{align}
    Inequality (\ref{explain eval}) holds because the left-hand side is the evaluation of the right-hand side at $y=0$. Hence,
    \begin{align*}
        \min_{w_{P}\in W_{P}} \E_{(x,c)\sim P}(\hat{c}_\theta (x)^\top  w_{P}(x))&\leq \E_{(x,c)\sim P}(\hat{c}_\theta (x)^\top  w(c(x))) + y^\star_{P}(\theta) (\E_{(x,c)\sim P}(c(x)^\top w(c(x)))-\beta)\\
        &=\E_{(x,c)\sim P}(\hat{c}_\theta (x)^\top  w(c(x))) + y^\star_{P}(\theta) (\beta_{\text{min},P}-\beta).
    \end{align*}
    This yields
    \begin{align*}
        y^\star_{P}(\theta)\leq \frac{\E_{(x,c)\sim P}(\hat{c}_\theta (x)^\top  w(c(x))) - \E_{(x,c)\sim P}(\hat{c}_\theta (x)^\top  w(\hat{c}_\theta(x)))}{\beta - \beta_{\text{min},P}}\leq \frac{2DB_WB_{\Phi}}{\beta - \beta_{\text{min},P}}.
    \end{align*}
    Hence, we get the following upper bound for the random variables we are working with
    \begin{align*}
        \abs{m(x,\theta)}\leq 2DB_WB_{\Phi} + KB_W \frac{2DB_WB_{\Phi}}{\beta - \beta_{\text{min},P}}:=U
    \end{align*}
    we can apply Hoeffding's inequality, and say that 
    \begin{align*}
        \mathbb P\left(\Tilde{\ell}_{P}^{\beta}(\theta)-\Tilde{\ell}_{P_n}^{\beta}(\theta)  >\varepsilon\right)\leq 2e^{-\frac{n\varepsilon^2}{2U^2}}.
    \end{align*}
    This yields
    \begin{align*}
        \mathbb P\left(\sup_{B\in \mathcal B}\left( \Tilde{\ell}_{P}^{\beta}(\theta_{B}) - \Tilde{\ell}_{P_n}^{\beta}(\theta_{B})\right)>\varepsilon\right)\leq \abs{\mathcal B}2e^{-\frac{n\varepsilon^2}{2U^2}}.
    \end{align*}
    By denoting $\delta := \abs{\mathcal B}2e^{-\frac{n\varepsilon^2}{2U^2}}$, we can say that with probability at least $1-\delta$, we have 
    \begin{align*}
        \forall B \in \mathcal B,\; \Tilde{\ell}_{P}^{\beta}(\theta_{B}) -\Tilde{\ell}_{P_n}^{\beta}(\theta_{B}) \leq U\sqrt{2}\sqrt{\frac{\log \abs{\mathcal B} + \log \frac{2}{\delta} }{n}},
    \end{align*}
    which gives 
    \begin{align*}
        \forall B \in \mathcal B,\; \ell_{P}^{\beta}(\theta_{B}) -\ell_{P_n}^{\beta}(\theta_{B})\leq U\sqrt{2}\sqrt{\frac{\log \abs{\mathcal B} + \log \frac{2}{\delta} }{n}}\leq U\sqrt{2}\sqrt{\frac{d \log \left(1+\frac{2D}{\gamma}\right) + \log \frac{2}{\delta} }{n}}.
    \end{align*}
    Now we have for every $\theta \in B_{\norm{.}}(0,D)$, denoting $B\in \mathcal B$ such that $\theta \in B$,
    \begin{align}
        \ell_{P}^{\beta}(\theta) -\ell_{P_n}^{\beta}(\theta) &\leq \ell_{P}^{\beta}(\theta_B) -\ell_{P_n}^{\beta}(\theta_B) + \abs{\ell_{P_n}^{\beta}(\theta) -\ell_{P}^{\beta}(\theta)  - \left(\ell_{P_n}^{\beta}(\theta_B) -\ell_{P}^{\beta}(\theta_B) \right)}\\
        &\leq  U\sqrt{2}\sqrt{\frac{d \log \left(1+\frac{2D}{\gamma}\right) + \log \frac{2}{\delta} }{n}} + 4B_WB_{\Phi}\norm{\theta - \theta_B}\\
        &\leq U\sqrt{2}\sqrt{\frac{d \log \left(1+\frac{2D}{\gamma}\right) + \log \frac{2}{\delta} }{n}} + 4B_WB_{\Phi}\gamma. \label{applying sensitivity}
    \end{align}
    Inequality (\ref{applying sensitivity}) is resulting from the Lipschitz inequality (\ref{gridineq}). Taking $\gamma = O\left(\sqrt{\frac{\log \frac{1}{\delta}}{n}}\right)$, we get 
    \begin{align*}
         \ell_{P}^{\beta}(\theta) -\ell_{P_n}^{\beta}(\theta)\leq O\left(\left( U+4B_WB_{\Phi}\right)\sqrt{\frac{\log \frac{1}{\delta}}{n}}\right).
    \end{align*}
    Since we have $\ell_{P_n}^{\beta}(\theta^\star)\leq \varepsilon$, we can replace in the inequality above $\theta$ by $\theta^\star$ and get
    \begin{align*}
         \ell_{P}^{\beta}(\theta^\star) \leq  \ell_{P_n}^{\beta}(\theta^\star) +O\left(\left( U+4B_WB_{\Phi}\right)\sqrt{\frac{\log \frac{1}{\delta}}{n}}\right)\leq \varepsilon +O\left(\frac{1}{\beta - \beta_{\text{min},P}}\sqrt{\frac{\log \frac{1}{\delta}}{n}}\right).
    \end{align*}
    \end{Proofe}
\subsection{Proof of Lemma \ref{positivity}}
\begin{Proofe}
    Since the only difference between the right and left optimization problems in $\ell_P^\beta$ is an additional constraint added to the right problem, it is clear that it will yield a higher value than the right one for any $\theta \in \R^m$, and consequently we have that indeed $\ell_P^\beta$ is positive.
\end{Proofe}
\subsection{Proof of Proposition \ref{regularity}}
\begin{Proofe}
For any $\theta \in \R^m\setminus \{0\},$ the uniqueness of the solution to the linear program $\min_{w\in W}\hat{c}_\theta (x)^\top w$ is a consequence of the condition $\prod_{v\in \mathcal V}^{}\hat{c}_\theta(x)^\top v \neq 0$. Since the set $\mathcal Z$ has a zero Lebesgue measure and that $x$ has a continuous distribution, it follows that 
$$\mathbb P\left(\prod_{v\in \mathcal V}^{}\hat{c}_\theta(x)^\top v = 0\right)=0,$$
i.e. when $\theta \neq 0$, the linear program $\min_{w\in W}\hat{c}_\theta (x)^\top w$ has a unique solution almost surely.
\end{Proofe}
\subsection{Proof of Proposition \ref{DC structure}} \label{DC property proof}
\begin{Proofe}
    We denote for every $w\in W$ and $\theta \in \R^m$, $f(\theta;w)=-\hat{c}_\theta(x)^\top w $. We want to prove that for every $w\in W$, $\theta \longrightarrow f(\theta;w)$ is weakly convex. More precisely, we want to prove that there exists $\alpha\geq 0$ such that for every $w\in W$, $\theta \longmapsto f(\theta;w)+\alpha \norm{\theta}^2$ is convex. We have for every $\theta_1,\theta_2 \in \R^m$, using assumptions \ref{smoothness} and \ref{boundedness},
    \begin{align*}
        \norm{\nabla f(\theta_1;w)-\nabla f(\theta_2;w)}&=\norm{\left(\nabla \hat{c}_{\theta_1} (x)-\nabla \hat{c}_{\theta_2} (x)\right)^\top w}\\
        &\leq B_W\norm{\nabla \hat{c}_{\theta_1} (x)-\nabla \hat{c}_{\theta_2} (x)}\\
        &\leq B_WB_L\norm{\theta_1-\theta_2}.
    \end{align*}
    We denote for every $\theta \in \R^m$ and $w\in W$, $h(\theta;w)=f(\theta;w)+2B_WB_L \norm{\theta}^2$. We have for every $\theta_1,\theta_2\in \R^m$,
    \begin{align*}
        (\nabla h(\theta_1;w)-\nabla h(\theta_2;w))^\top (\theta_1-\theta_2)&=(\nabla f (\theta_1;w)-\nabla f (\theta_2;w))^\top (\theta_1-\theta_2) + 2B_WB_L\norm{x-y}^2\\
        &\geq -\norm{\nabla f (\theta_1;w)-\nabla f (\theta_2;w)}\norm{ \theta_1-\theta_2}+2B_WB_L \norm{\theta_1-\theta_2}^2 \\
        &\geq -B_WB_L\norm{\theta_1-\theta_2}^2+2B_WB_L \norm{\theta_1-\theta_2}^2 \\
        &=B_WB_L\norm{\theta_1-\theta_2}^2\geq 0.
    \end{align*}
    Hence, for every $w\in W$, $\theta \longmapsto h(\theta;w)$ is convex. This implies that for any measurable mapping $w:\R^k\longrightarrow W$, $\theta \longmapsto \E_{(x,c)\sim P_n}\left(h(\theta;w(x))\right)$ is convex. Consequently, $\theta \longmapsto \max_{w_{P_n}\in W_{P_n}}\E_{(x,c)\sim P_n}(h(\theta;w(x)))$ and $\theta \longmapsto \max_{w_{P_n}\in \overline W_{P_n}^\beta}\E_{(x,c)\sim P_n}(h(\theta;w(x)))$ are both maximums of a family of convex functions, and hence are also convex functions. Finally we have for every $\theta \in \R^m$,
    \begin{align*}
        \max_{w_{P_n}\in W_{P_n}}\E_{(x,c)\sim P_n}(h(\theta;w(x))) = g_{P_n}(\theta) + 2B_L\norm{\theta}^2 \\ \max_{w_{P_n}\in \overline W_{P_n}^\beta}\E_{(x,c)\sim P_n}(h(\theta;w(x))) = \overline g_{P_n}^\beta(\theta) + 2B_L\norm{\theta}^2.
    \end{align*}
    Finally, we can write for every $\theta \in \R^m$,
    \begin{align*}
        \ell_{P_n}^\beta(\theta)=g_{P_n}(\theta) + 2B_L\norm{\theta}^2 - \left( \overline g_{P_n}^\beta(\theta) + 2B_L\norm{\theta}^2\right) 
    \end{align*}
    which is indeed a difference of convex functions.
\end{Proofe}

\subsection{Proof of Theorem \ref{sensitivity}} \label{sensitivity proof}
\begin{Proofe}
Let $\theta \in \R^m$ and $t>0$. Assumption \ref{marginal hypothesis} gives 
\begin{align*}
    \mathbb P(0<\Delta_\theta(x)\leq \norm{\theta}t)\leq \left(\frac{\gamma t}{B}\right)^\alpha.
\end{align*}
where 
\begin{align*}
    \Delta_\theta(x)=\begin{cases}
             \min_{w \in W^\angle \setminus W^\star (\hat{c}_\theta(x))}\hat{c}_\theta(x)^\top w- \min_{w \in W^\angle}\hat{c}_\theta(x)^\top w & \text{ if }W^\star (\hat{c}_\theta(x))\neq W\\
            0  & \text{ else. }\end{cases}
\end{align*}
We have for $\theta \neq 0$, denoting $W^\angle =\{w_1,\dots,w_N\}$,  $N\in \mathbb N$,
\begin{align*}
    &\Delta_\theta (x) > \norm{\theta}t \Longrightarrow \forall w_\angle \in W^\angle \setminus \{w(\hat{c}_\theta)\},\;  \hat{c}_\theta(x)^\top w_\angle -\hat{c}_\theta(x)^\top w(\hat{c}_\theta(x))>\norm{\theta} t\\
    &\Longrightarrow \forall w_\angle \in W^\angle \setminus \{w(\hat{c}_\theta)\},\;  \hat{c}_\theta(x)^\top w_\angle - \hat{c}_\theta(x)^\top w(\hat{c}_\theta(x)) >\norm{\theta}\frac{\norm{w(\hat{c}_\theta (x)) - w_\angle}}{2B_W} t
    \\&\Longrightarrow \forall \lambda_1,\dots,\lambda_n \geq 0\; \text{s.t. }\lambda_1+\dots+\lambda_n=1,\\& \sum_{i=1}^{N}\lambda_i \left(\hat{c}_\theta(x)^\top w_i - \hat{c}_\theta(x)^\top w(\hat{c}_\theta(x)) \right) >\frac{1}{2B_W}\norm{\theta}t\sum_{i=1}^{N}\lambda_i\norm{w(\hat{c}_\theta (x)) - w_i}
    \\&\Longrightarrow \forall \lambda_1,\dots,\lambda_n \geq 0\; \text{s.t. }\lambda_1+\dots+\lambda_n=1,\\&  \hat{c}_\theta(x)^\top \sum_{i=1}^{N}\lambda_iw_i-\hat{c}_\theta(x)^\top w(\hat{c}_\theta(x)) >\frac{1}{2B_W}\norm{\theta}t\norm{w(\hat{c}_\theta (x)) - \sum_{i=1}^{N}\lambda_iw_i}\\
    &\Longrightarrow \forall w\in W,\; \hat{c}_\theta(x)^\top w-\hat{c}_\theta(x)^\top w(\hat{c}_\theta(x)) >\frac{1}{2B_W}\norm{\theta}\norm{w(\hat{c}_\theta (x)) - w}t.
\end{align*}
Hence, we have 
\begin{align}
    \mathbb P \left(\forall w\in W,\; \hat{c}_\theta(x)^\top w - \hat{c}_\theta(x)^\top w(\hat{c}_\theta(x)) >\frac{1}{2B_W}\norm{\theta}\norm{w(\hat{c}_\theta (x)) - w}t\right)\geq \mathbb P \left(\Delta_\theta (x)>\norm{\theta}t\right)>1-\left(\frac{\gamma t}{B}\right)^\alpha \label{probability bound}
\end{align}
We denote for $w\in W$, $t\geq 0$, and $\theta \in \R^m \setminus \{0\}$, 
\begin{align*}
    G(x,\theta,w,t)= \hat{c}_\theta(x)^\top w - \hat{c}_\theta(x)^\top w(\hat{c}_\theta(x))  -\frac{1}{2B_W}\norm{\theta}\norm{w(\hat{c}_\theta (x)) - w}t
\end{align*}
and the event $A(x,\theta,t)$ defined by 
\begin{align*}
    A(x,\theta,t)=\left\{\forall w\in W,\; G(x,\theta,w(x),t)>0\right\}.
\end{align*}
Hence, we have for every mapping $w:\R^k\longrightarrow w$ and $\theta \in \R^m\setminus\{0\}$,
\begin{align*}
    \E_{(x,c)\sim P}(G(x,\theta,w(x),t))&=\E_{(x,c)\sim P}(G(x,\theta,w(x),t)1_{A(x,\theta,t)}) + \E_{(x,c)\sim P}(G(x,\theta,w(x),t)1_{A(x,\theta,t)^c})\\
    &\geq \E_{(x,c)\sim P}(G(x,\theta,w(x),t)1_{A(x,\theta,t)^c}).
\end{align*}
Furthermore, we have for every $w\in W$
\begin{align*}
    \abs{G(x,\theta,w,t)}& = \abs{ \hat{c}_\theta(x)^\top w - \hat{c}_\theta(x)^\top w(\hat{c}_\theta(x))  -\frac{1}{2B_W}\norm{\theta}\norm{w(\hat{c}_\theta (x)) - w}t}
    \\&\leq B_{\Phi} B_W\norm{\theta} + B_\Phi B_W \norm{\theta} + \norm{\theta}t=\norm{\theta} \left(2B_W B_\Phi +t\right).
\end{align*}
Notice that (\ref{probability bound}) can be rewritten as $\mathbb P \left(A(x,\theta,t)\right)>1-\left(\frac{\gamma t}{B}\right)^\alpha$. Hence, we have
\begin{align*}
    \E_{(x,c)\sim P}(G(x,\theta,w(x),t))&\geq - (2B_W B_\Phi +t)\norm{\theta} \mathbb P (A(x,\theta,t)^c)
    \\&\geq - (2B_W B_\Phi +t)\norm{\theta} \left(\frac{\gamma t}{B}\right)^\alpha.
\end{align*}
In conclusion, we have for every $t\geq 0$ and every $\theta \in \R^m$, and mapping $w:\R^k \longrightarrow W$,
\begin{align*}
    \E_{(x,c)\sim P}\left(\hat{c}_\theta(x)^\top w(x) - \hat{c}_\theta(x)^\top w(\hat{c}_\theta(x)) \right) \geq \frac{1}{2B_W}\norm{\theta}\E_{(x,c)\sim P}\left(\norm{w(\hat{c}_\theta (x)) - w(x)}\right)t- (2B_W B_\Phi +t)\norm{\theta} \left(\frac{\gamma t}{B}\right)^\alpha.
\end{align*}
the inequality above is trivially verified when $\theta =0$. Also, when taking $w\in \arg\min_{w\in \overline W_P^\beta}\E_{(x,c)\sim P}\left(\hat{c}_\theta(x)^\top w(x)\right)$, we get
\begin{align*}
    \ell_P^\beta(\theta)&\geq \frac{1}{2B_W}\norm{\theta}\E_{(x,c)\sim P}\left(\norm{w(\hat{c}_\theta (x)) - w(x)}\right)t- (2B_W B_\Phi +t)\norm{\theta} \left(\frac{\gamma t}{B}\right)^\alpha\\
    &\geq \frac{\norm{\theta}t}{2KB_W}(\ell_P(\theta) - \beta ) - (2B_W B_\Phi +t)\norm{\theta} \left(\frac{\gamma t}{B}\right)^\alpha
\end{align*}
Let $\varepsilon>0$. Assume that $\theta$ is bounded away from $0$. We have
\begin{align*}
    \ell_{P}(\theta) - \beta\leq \frac{2KB_W}{t}\left(\frac{\ell_P^\beta(\theta)}{\norm{\theta}} + (2B_W B_\Phi +t) \left(\frac{\gamma t}{B}\right)^\alpha \right)
    \leq  O\left(\frac{1}{t}\left(\ell_P^\beta \left(\theta\right) + t^\alpha\right)\right)
\end{align*}
 To optimize the right hand side and ignoring the constant, we take $t=\ell_P^\beta \left(\theta\right) ^{\frac{1}{1+\alpha}}$. This gives
 \begin{align*}
     \ell_P(\theta)\leq \beta + O\left(\ell_P^\beta \left(\theta\right) ^{1-\frac{1}{1+\alpha}}\right).
 \end{align*}
\end{Proofe}
\subsection{Proof of Theorem \ref{subdifferentialsurrogate}} \label{subdifferentialsurrogate proof}
\begin{Proofe}{}
\begin{enumerate}
    \item We denote for all $(\theta,v)\in \R^m\times V_{P}$, $\phi(\theta,v)=\theta^\top v$, $v_\beta$ and for all $(\theta,v_\beta)\in \R^m\times V_{P}$, $\phi(\theta,v_\beta)=\theta^\top v_\beta$,  In this case, for all $\theta \in \R^m$, the loss writes as
\begin{align*}
    \ell_{P}^\beta(\theta)=g_{P}(\theta)-\overline{g}_{P}^\beta(\theta)=\min_{v_\beta\in V_{P}^\beta}\phi(\theta,v_\beta) - \min_{v\in V_{P}}\phi(\theta,v).
\end{align*}
Let $\mathcal V_{P}^\beta(\theta)$ and $\mathcal V_{P}(\theta)$ be the set of minimizers of respectively $v_\beta\longmapsto \phi(\theta,v_\beta)$ over $V_{P}^\beta$ and $v\longmapsto \phi(\theta,v)$ over $V_{P}$. Given that $V_{P}$ and $V_{P}^\beta$ are compact sets, and $\theta \longmapsto \phi(\theta,v)$ for all $v\in V_{P}\cup V_{P}^\beta$ is differentiable, and $\phi$ is continuous, we can say  by Danskin's theorem that for all $\theta \in \R^m$,
\begin{align*}
    \partial g_{P}(\theta)&=\text{conv}\left\{\frac{\partial \phi(\theta,v)}{\partial \theta },\; v\in \mathcal V_{P}(\theta)\right\}=\text{conv} \;\mathcal V_{P}(\theta)\underset{(*)}{=}\mathcal V_{P}(\theta),\\
    \partial \overline{g}^\beta_{P}(\theta)&=\text{conv}\left\{\frac{\partial \phi(\theta,v_\beta)}{\partial \theta },\; v_\beta\in \mathcal V_{P}^\beta(\theta)\right\}=\text{conv} \;\mathcal V_{P}^\beta(\theta)\underset{(**)}{=}\mathcal V_{P}^\beta(\theta).
\end{align*}
The two inequalities $(*)$ and $(**)$ are due to the fact that $\mathcal V_{P}(\theta)$ and $\mathcal V_{P}^\beta(\theta)$ are convex sets. Furthermore, we have for all $\theta \in \R^m$,
\begin{align*}
    \ell_{P}^\beta(\theta)=\min_{v_\beta\in V_{P}^\beta}\phi(\theta,v_\beta) - \min_{v\in V_{P}}\phi(\theta,v) = \theta^\top  v^\star_\beta - \theta ^\top v^\star=\theta^\top \underbrace{\left(v^\star_\beta-v^\star\right)}_{:=v(\theta)},
\end{align*}
where $(v^\star_\beta,v^\star)\in \mathcal V_{P}^\beta(\theta)\times \mathcal V_{P}(\theta)$. All of the above clearly yields $$v(\theta)\in \mathcal V_{P}^\beta(\theta)-\mathcal V_{P}(\theta)=\partial \overline{g}_{P}^\beta(\theta)-\partial g_{P}(\theta),$$
which is the result we were seeking to prove.
\item From the above, it is clear that for a given $\theta \in \R^m$ if $v(\theta)=0$, then $\ell_{P}^\beta(\theta)=0$. Furthermore, using the very first definition of $\ell^\beta_{P}$, since $\overline{W}_P^\beta \subset W_P$, we have for all $\theta \in \R^m$,
\begin{align*}
    \ell_{P}^\beta(\theta)&=\min_{{w}_{P}\in \overline{W}_{P}^{\beta}}\mathbb E_{(c,x)\sim P}\left((\hat{c}_{\theta}(x)^\top w_{P}(x))\right)-\min_{{w}_{P}\in W_{P}}\mathbb E_{(c,x)\sim P}\left(\hat{c}_{\theta}(x)^\top w_{P}(x)\right)\geq 0.
\end{align*}
In conclusion, $0$ is a lower bound of $\ell_{P}^\beta$, and if for a given $\theta$, $v(\theta)=0$, then $\ell_{P}^\beta(\theta)=0$, i.e. $\theta$ is a minimizer of $\ell_{P}^\beta$.
\item $r_{P}^\beta$ is the difference of the Moreau envelope of two convex functions, hence it is the difference between two differentiable convex functions. this yields that $r_{P}^\beta$ is differentiable. Furthermore, if $\lambda$ is a stationary point of $r_{P}^\beta$ then using proposition 1 from \cite{sun2021algorithms}, $\theta^\beta_{P}(\lambda)$ is indeed a stationary point of $\ell_{P}^\beta$. Furthermore, we have for every $\varepsilon\geq 0$, We consider $\lambda \in \R^m$ an $\varepsilon$-stationary point of $r_{P_n}^{\beta}$, i.e. $\norm{\nabla r_{P_n}^{\beta} (\lambda)}\leq \varepsilon$. Using Danskin's theorem, we get
    \begin{align*}
        \nabla r_{P_n}^{\beta} (\lambda)&=\nabla M_{P_n}(\lambda) - \nabla \overline{M}_{P_n}^\beta (\lambda) \\& =\lambda - \theta_{P_n}(\lambda) - (\lambda - \overline{\theta}_{P_n}^{\beta}(\lambda))
        \\&= \overline{\theta}_{P_n}^{\beta}(\lambda)- \theta_{P_n}(\lambda),
    \end{align*}
    Where $\overline{\theta}_{P_n}^{\beta}(\lambda)=\arg \min_{\theta \in \R^m}g_{P_n}^{\beta}(\theta) +\frac{1}{2}\norm{\lambda - \theta}^2$ and $\theta_{P_n}(\lambda)=\arg \min_{\theta \in \R^m}g_{P_n}(\theta) +\frac{1}{2}\norm{\lambda - \theta}^2$. By looking at the proof of Lemma \ref{transferoptimality}, we can easily see that we  also have $\theta_{P_n}(\lambda)=\lambda +v_{P_n}(\theta_{P_n}(\lambda))$ where $v_{P_n}(\theta_{P_n}(\lambda))\in \arg \min_{v\in V_{P_n}} \theta_{P_n}(\lambda)^\top v=\partial g_{P_n}(\theta_{P_n}(\lambda))$ and $\overline{\theta}_{P_n}^{\beta}(\lambda)=\lambda +v_{P_n}^{\beta}(\overline{\theta}_{P_n}^{\beta}(\lambda))$ where $v_{P_n}^{\beta}(\overline{\theta}_{P_n}^{\beta}(\lambda))\in \arg \min_{v\in \overline{V}^\beta_{P_n}} \overline{\theta}_{P_n}^{\beta}(\lambda)^\top v = \partial g_{P_n}^{\beta}(\overline{\theta}_{P_n}^{\beta}(\lambda))$. Let $v_{P_n}(\overline{\theta}_{P_n}^{\beta}(\lambda))\in \arg \min_{v\in V_{P_n}} \overline{\theta}_{P_n}^{\beta}(\lambda)^\top v=\partial g_{P_n}(\overline{\theta}_{P_n}^{\beta}(\lambda))$. We have 
    \begin{align}
        \ell_{P_n}^\beta(\overline{\theta}_{P_n}^{\beta}(\lambda))&=\overline{\theta}_{P_n}^{\beta}(\lambda)^\top v_{P_n}^\beta(\overline{\theta}_{P_n}^{\beta}(\lambda))-\overline{\theta}_{P_n}^{\beta}(\lambda)v_{P_n}(\overline{\theta}_{P_n}^{\beta}(\lambda))\\
        &\leq \abs{\overline{\theta}_{P_n}^{\beta}(\lambda)^\top v_{P_n}^\beta(\overline{\theta}_{P_n}^{\beta}(\lambda))-\overline{\theta}_{P_n}^{\beta}(\lambda)^\top v_{P_n}(\theta_{P_n}(\lambda))}\\&+\abs{\theta_{P_n}(\lambda)^\top v_{P_n}(\theta_{P_n}(\lambda))-\overline{\theta}_{P_n}^{\beta}(\lambda)^\top v_{P_n}(\overline{\theta}_{P_n}^{\beta}(\lambda))}\\&+\abs{\theta_{P_n}(\lambda)^\top v_{P_n}(\theta_{P_n}(\lambda))-\overline{\theta}_{P_n}^{\beta}(\lambda)^\top v_{P_n}(\theta_{P_n}(\lambda))}\\&\leq  \norm{\overline{\theta}_{P_n}^{\beta}(\lambda)}\norm{v_{P_n}^\beta(\overline{\theta}_{P_n}^{\beta}(\lambda))-v_{P_n}(\theta_{P_n}(\lambda))}+\abs{\kappa(\theta_{P_n}^\beta(\lambda)) - \kappa(\theta_{P_n}(\lambda))}\\&+\norm{v_{P_n}(\theta_{P_n}(\lambda))}\norm{\theta_{P_n}(\lambda)-\overline{\theta}_{P_n}^{\beta}(\lambda)}\\
        &\leq (2B_WB_\Phi +5B_WB_{\Phi})\varepsilon + \abs{\kappa(\theta_{P_n}^\beta(\lambda)) - \kappa(\theta_{P_n}(\lambda))}
        \\&\leq (3B_WB_\Phi +5B_WB_{\Phi})\varepsilon =8B_WB_\Phi \varepsilon\label{lasteq}
    \end{align}
    Here, $\kappa$ is defined for any $\theta \in \R^m$ as $\kappa(\theta)=\min_{v\in V_{P_n}}\theta^\top v$. The last equality \ref{lasteq} is due to the fact that the subgradient of $\kappa$ at some $\theta \in \R^m$ is in $V_{P_n}$ (which can be proven thanks to Danskin's theorem) and is hence bounded by $B_WB_\Phi$. We can see that inequality \ref{lasteq} is indeed the inequality we were seeking to obtain. Taking $\varepsilon= 0$ gives that stationarity for $r_{P_n}^\beta$ implies global optimality for $\ell_{P_n}^\beta$.
\end{enumerate}

\end{Proofe}
\subsection{Proof of Lemma \ref{transferoptimality}} \label{transferoptimality proof}
\begin{Proofe}{}
The main idea of the proof is to switch min and max in the definitions of $M_{P}^{\beta}$ and $\overline{M}_{P}^{\beta}$.
\begin{enumerate}
    \item We have for all $\lambda \in \mathbb R^m$,
    \begin{align}
        M_{P}(\lambda)&=\min_{\theta \in \mathbb R^m}g_{P}^\beta(\theta) +\frac{1}{2} \norm{\lambda -\theta}^2\\
        &=\min_{\theta \in \mathbb R^m}\left(-\min_{w_{P} \in W_{P}}\E (\hat{c}_\theta (x)^\top w_{P}(x))\right) +\frac{1}{2} \norm{\lambda -\theta}^2\\
        &=\min_{\theta \in \mathbb R^m}\max_{w_{P} \in W_{P}}-\E (\hat{c}_\theta (x)^\top w_{P}(x))+\frac{1}{2} \norm{\lambda -\theta}^2\\
        &=\min_{\theta \in \mathbb R^m}\max_{w_{P} \in W_{P}}-\theta^\top \underbrace{\E (\Phi(x) ^\top w_{P}(x))}_{v}+\frac{1}{2} \norm{\lambda -\theta}^2\\
        &=\min_{\theta \in \mathbb R^m}\max_{v \in V_{P}}-\theta^\top v+\frac{1}{2} \norm{\lambda -\theta}^2\label{minmaxbefore}\\
        &=\max_{v \in V_{P}}\min_{\theta \in \mathbb R^m}-\theta^\top v+\frac{1}{2} \norm{\lambda -\theta}^2 \label{invminmax} 
    \end{align}
    The equality between \ref{invminmax} and \ref{minmaxbefore} holds because of Sion's minimax theorem: since $\R^m$ and $V_{P}$ are convex sets, $v\longmapsto -\theta^\top  v$ is upper semi-continuous and concave for any $\theta \in \R^m$, and $\theta\longmapsto -\theta^\top  v +\frac{1}{2}\norm{\lambda -\theta}^2$ is lower semicontinuous and convex for any $v\in V_{P}$, we can switch min and max. In \ref{invminmax}, the minimum is reached when the gradient with respect to $\theta$ is zero, i.e. $-v+\theta - \lambda =0$, which gives $\theta = v+\lambda$. In this case, we get
    \begin{align*}
        M_{P}(\lambda)&=\max_{v \in V_{P}}-(v+\lambda)^\top v+\frac{1}{2} \norm{v}^2\\
        &=\max_{v \in V_{P}}-\lambda^\top v-\frac{1}{2} \norm{v}^2\\
        &=\frac{1}{2}\norm{\lambda}^2 + \max_{v \in V_{P}}-\frac{1}{2}\norm{\lambda + v}^2\\
        &=\frac{1}{2}\norm{\lambda}^2 - \min_{v \in V_{P}}\frac{1}{2}\norm{\lambda + v}^2.
    \end{align*}
    The calculations above also give us immediately that $\theta_{P}(\lambda)=\lambda + v$ where $v=\arg \min_{v \in V_{P}}\frac{1}{2}\norm{\lambda + v}^2$.
    \item The proof of the second property is almost identical to the previous one. It suffices to replace $V_{P}$ by $\overline{V}_{P}^\beta$ and $W_{P}$ by $\overline{W}_{P}^\beta$.
\end{enumerate}
\end{Proofe}
\subsection{Proof of Proposition \ref{representation}} \label{representation proof}
This proposition immediately follows from Lemma \ref{transferoptimality}.
\subsection{Proof of Theorem \ref{move}}
\begin{Proofe}{}
    We will sequentially prove the two points. We assume that $\lambda \in \R^m$ is a stationary point of $f_{P_n}^\beta$. In this case, we have 
    \begin{align*}
        0=\nabla f_{P_n}^\beta(\lambda) = \frac{\lambda - F_{\overline V_{P_n}^\beta}(\lambda)}{\norm{\lambda - F_{\overline V_{P_n}^\beta}(\lambda)}^2} -\frac{\lambda - F_{V_{P_n}}(\lambda)}{\norm{\lambda - F_{V_{P_n}}(\lambda)}^2}.
    \end{align*}
    This yields
    \begin{align*}
        \frac{\lambda - F_{\overline V_{P_n}^\beta}(\lambda)}{\norm{\lambda - F_{\overline V_{P_n}^\beta}(\lambda)}^2}=\frac{\lambda - F_{V_{P_n}}(\lambda)}{\norm{\lambda - F_{V_{P_n}}(\lambda)}^2}.
    \end{align*}
    Taking the norm in both sides, we get $\norm{\lambda - F_{\overline V_{P_n}^\beta}(\lambda)}=\norm{\lambda - F_{V_{P_n}}(\lambda)}$, and finally, replugging this in the equality above, we directly get $F_{\overline V_{P_n}^\beta}(\lambda)=F_{V_{P_n}}(\lambda)$, i.e. $\nabla r_{P_n}^\beta(\lambda)=0$. Let us now prove the second result. We denote $D_{n}^\beta=\lambda - F_{\overline V_{P_n}^\beta}(\lambda)$ and $D_n=\lambda - F_{V_{P_n}}(\lambda)$. We assume now that $\norm{f_{P_n}^\beta(\lambda)}\leq \varepsilon$. This inequality gives us two inequalities
    \begin{align*}
        \abs{\frac{1}{\norm{D_n}}-\frac{1}{\norm{D_n^\beta}}}\leq \varepsilon \text{ and } \norm{\frac{D_n}{\norm{D_n}^2}-\frac{D_n^\beta}{\norm{D_n^\beta}^2}}\leq \varepsilon.
    \end{align*}
    Hence, we have 
    \begin{align*}
        \norm{\nabla r_{P_n}^\beta}&=\norm{F_{\overline V_{P_n}^\beta}(\lambda) - F_{V_{P_n}}(\lambda)}\\
        &=\norm{D_n^\beta - D_n}\\
        &=\norm{D_n}^2 \norm{\frac{D_n}{\norm{D_n}^2}-\frac{D_n^\beta}{\norm{D_n}^2}}\\
        &\leq \norm{D_n}^2 \left(\norm{\frac{D_n}{\norm{D_n}^2}-\frac{D_n^\beta}{\norm{D_n^\beta}^2}} + \norm{D_n^\beta}\abs{\frac{1}{\norm{D_n}^2} - \frac{1}{\norm{D_n^\beta}^2}}\right)\\
        &\leq 9B_W^2B_\Phi^2 \varepsilon + 9B_W^2B_\Phi^2\left(1+\frac{\norm{D_n}}{\norm{D_n^\beta}}\right)\abs{\frac{1}{\norm{D_n}}-\frac{1}{\norm{D_n^\beta}}}\\
        &\leq 27 B_W^2B_\Phi^2 \varepsilon.
    \end{align*}

    Let $\lambda \in \R^m$ such that $\norm{\nabla r_{P_n}^\beta (\lambda)}\leq \varepsilon$, i.e. \begin{align}\label{epsstat}
        \norm{F_{\overline V_{P_n}^\beta}(\lambda)-F_{V_{P_n}}(\lambda)}\leq \varepsilon.
    \end{align} We have 
     \begin{align}
         \norm{\Tilde{\lambda} - F_{\overline V_{P_n}^\beta}(\Tilde{\lambda})}&\leq \norm{\Tilde{\lambda} - F_{\overline V_{P_n}^\beta}(\lambda)}\\
         &\leq \norm{\Tilde{\lambda} - F_{V_{P_n}}(\lambda)} + \norm{F_{V_{P_n}}(\lambda) - F_{\overline V_{P_n}^\beta}(\lambda)}\\
         &\leq \norm{\Tilde{\lambda} - F_{V_{P_n}}(\Tilde{\lambda})} + \varepsilon.\label{epsilonineq}
     \end{align}
The last inequality is due to inequality \ref{epsstat} and $F_{V_{P_n}}(\Tilde{\lambda})=F_{\overline V_{P_n}^\beta}(\lambda)$. We finish by proving this last equality. We have for every $v\in -V_{P_n}$,
\begin{align}
    \norm{\Tilde{\lambda} - v}^2&=\norm{\lambda + r \frac{\lambda-F_{V_{P_n}}(\lambda)}{\|\lambda-F_{V_{P_n}}(\lambda)\|} - v}^2\\
    &=\norm{\lambda - v}^2 + r^2 +2r\ps{\lambda - v }{\frac{\lambda-F_{V_{P_n}}(\lambda)}{\|\lambda-F_{V_{P_n}}(\lambda)\|}}\\ &\geq 
    \norm{\lambda - v}^2 + r^2. \label{FINALLHS}
\end{align}
Notice that when the left-hand side in \ref{FINALLHS} is minimized (i.e. $v=F_{V_{P_n}}(\lambda)$, the right-hand side takes the same value. Hence, $v=F_{V_{P_n}}(\lambda)$ is also a minimizer of the left-hand side, i.e. $F_{V_{P_n}}(\Tilde{\lambda})=F_{\overline V_{P_n}^\beta}(\lambda)$. Finally, inequality \ref{epsilonineq} yields
\begin{align*}
    \norm{\Tilde{\lambda} - F_{\overline V_{P_n}^\beta}(\Tilde{\lambda})}-\norm{\Tilde{\lambda} - F_{V_{P_n}}(\Tilde{\lambda})} \leq  \varepsilon,
\end{align*}
Which gives 
\begin{align*}
    r_{P_n}^\beta (\Tilde{\lambda})&=\frac{1}{2}\norm{\Tilde{\lambda} - F_{\overline V_{P_n}^\beta}(\Tilde{\lambda})}^2-\frac{1}{2}\norm{\Tilde{\lambda} - F_{V_{P_n}}(\Tilde{\lambda})}^2\\
    &=\frac{1}{2}\left(\norm{\Tilde{\lambda} - F_{\overline V_{P_n}^\beta}(\Tilde{\lambda})}-\norm{\Tilde{\lambda} - F_{V_{P_n}}(\Tilde{\lambda})}\right) \left(\norm{\Tilde{\lambda} - F_{\overline V_{P_n}^\beta}(\Tilde{\lambda})}+\norm{\Tilde{\lambda} - F_{V_{P_n}}(\Tilde{\lambda})}\right)\\
    &\leq \frac{1}{2}\left(\norm{\Tilde{\lambda} - F_{\overline V_{P_n}^\beta}(\Tilde{\lambda})}+\norm{\Tilde{\lambda} - F_{V_{P_n}}(\Tilde{\lambda})}\right) \varepsilon
    \\&\leq \frac{1}{2}( 6B_WB_\Phi +2r) \varepsilon \leq 11B_WB_\Phi \varepsilon.
\end{align*}
Using this, we also prove the bound on the gradient of $r_{P_n}^\beta$. Using the properties of the projection $F_{V_{P_n}}(\Tilde{\lambda})$, we have
\begin{align*}
    \norm{\nabla r_{P_n}^\beta (\Tilde{\lambda})}&=\sqrt{\norm{F_{\overline V_{P_n}^\beta}(\Tilde{\lambda}) - F_{V_{P_n}}(\Tilde{\lambda})}^2}\\
    &=\sqrt{\norm{F_{\overline V_{P_n}^\beta}(\Tilde{\lambda}) - \Tilde{\lambda}}^2 - \norm{F_{V_{P_n}}(\Tilde{\lambda})-\Tilde{\lambda}}^2 -2\smash{\underbrace{\ps{F_{\overline V_{P_n}^\beta}(\Tilde{\lambda})- F_{V_{P_n}}(\Tilde{\lambda})}{F_{V_{P_n}}(\Tilde{\lambda})-\Tilde{\lambda}}}_{\geq 0}}}\\
    &\leq \sqrt{r_{P_n}^\beta (\tilde{\lambda})}\leq \sqrt{11B_WB_\Phi \varepsilon}.
\end{align*}
We finish by proving the bound on $\overline{\theta}_{P_n}^\beta (\Tilde{\lambda})$. Using Lemma \ref{representation}, we have $\overline{\theta}_{P_n}^\beta (\Tilde{\lambda})=\Tilde{\lambda}+F_{\overline V_{P_n}^\beta}(\Tilde{\lambda}) 
$, hence,
\begin{align*}
    B_WB_\Phi\leq \norm{\overline{\theta}_{P_n}^\beta (\Tilde{\lambda})}\leq 11B_WB_\Phi.
\end{align*}
\end{Proofe}

\end{document}